\documentclass[aap,seceqn,MSNbibl,citesort,dvips]{arximspdf}
\usepackage{mathrsfs}

\doi{10.1214/10-AAP709}
\volume{21}
\issue{2}
\pubyear{2011}
\firstpage{745}
\lastpage{799}

\makeatletter

\newtheorem{theorem}{Theorem}[section]
\newproclaim{defn}{Definition}[section]
\newtheorem{lem}{Lemma}[section]
\newtheorem{prop}{Proposition}[section]
\newtheorem{cor}{Corollary}[section]
\newproclaim{rem}{Remark}[section]

\newcommand{\R}{\mathbb{R}}
\newcommand{\Z}{\mathbb{Z}}
\newcommand{\N}{\mathbb{N}}
\newcommand{\M}{\mathbf{M}}
\newcommand{\D}{\mathbf{D}}
\newcommand{\E}{\mathbb{E}}
\newcommand{\pov}{\mathbf{d}}
\newcommand{\sko}{\varrho}
\newcommand{\init}{\mathscr I}
\newcommand{\ds}{\mathscr D}
\newcommand{\sol}{\Xi}
\newcommand{\solx}{\Upsilon}

\newcommand{\fr}[1]{\bar{#1}^r} 
\newcommand{\frm}[1]{\bar{#1}^{r,m}} 
\newcommand{\fro}[1]{\bar{#1}^{r,0}} 
\newcommand{\dr}[1]{\hat{#1}^r} 
\newcommand{\osc}[2]{\mathbf w_L({#1}(\cdot),#2)} 
\newcommand{\oscw}[2]{\mathbf w_T({#1},{#2})} 

\newcommand{\testfn}{\mathscr{V}}
\newcommand{\testset}{\mathcal{C}}

\newcommand{\ser}{\mathcal{Z}}
\newcommand{\buf}{\mathcal{Q}}
\newcommand{\dto}{\Rightarrow}

\newcommand{\gd}{\pi}
\newcommand{\K}{\mathbf K}

\makeatother

\begin{document}
\begin{frontmatter}

\title{Diffusion limits of limited processor sharing queues\thanksref{T1}}
\runtitle{Diffusion limits of LPS queues}

\thankstext{T1}{Supported in part by NSF Grants CMMI-0727400,
CMMI-0825840, CNS-0718701 and DMS-08-05979
and by two IBM Faculty awards and a VIDI grant from NWO.}

\begin{aug}
\author[A]{\fnms{Jiheng} \snm{Zhang}\corref{}\ead
[label=e1]{jiheng@ust.hk}},
\author[B]{\fnms{J. G.} \snm{Dai}\ead[label=e2]{dai@gatech.edu}} and
\author[C]{\fnms{Bert} \snm{Zwart}\thanksref{bert-aff}\ead
[label=e3]{bert.zwart@cwi.nl}}
\runauthor{J. Zhang, J. G. Dai and B. Zwart}
\affiliation{Hong Kong University of Science and Technology,
Georgia Institute of Technology and CWI}
\address[A]{J. Zhang\\
Department of Industrial Engineering\\
and Logistic Management\\
Hong Kong University of Science\\
and Technology\\
Clear Water Bay\\
Hong Kong\\
\printead{e1}}
\address[B]{J. G. Dai\\
H. Milton Stewart School\\
of Industrial and Systems Engineering\\
Georgia Institute of Technology\\
Atlanta, Georgia 30332\\
USA\\
\printead{e2}}
\address[C]{B. Zwart\\
CWI\\
P.O. Box 94079\\
1090 GB Amsterdam\\
The Netherlands\\
\printead{e3}}
\end{aug}

\thankstext{bert-aff}{This author is also affiliated with EURANDOM, VU
University
Amsterdam, and Georgia Institute of Technology.}

\received{\smonth{12} \syear{2007}}
\revised{\smonth{4} \syear{2010}}

%
\begin{abstract}
We consider a processor sharing queue where the number of jobs
served at any time is limited to $K$, with the excess jobs waiting
in a buffer. We use random counting measures on the positive axis
to model this system. The limit of this measure-valued process is
obtained under diffusion scaling and heavy traffic conditions. As a
consequence, the limit of the system size process is proved to be
a piece-wise reflected Brownian motion.
\end{abstract}

%

\end{frontmatter}

\section{Introduction}\label{sec:introduction}

This paper is concerned with developing a diffusion approximation for
a \textit{limited processor sharing} (LPS) queue, which consists of a
single server and an infinite capacity buffer. In such a system, the
server can serve up to $K\geq1$ jobs simultaneously, equally
distributing its attention to each of them. In other words, each job
in the server is processed at a rate that is the reciprocal of the
number of jobs in the server. An arriving job will immediately enter
the server and start receiving service if there are less than $K$ jobs
in the server when it arrives; otherwise it will wait in the buffer.
A job will leave the system immediately after the server has fulfilled
its service requirement. When the number of jobs in the server drops
from $K$ to $K-1$, the server will immediately admit the longest
waiting customer from the buffer, if there is one. We assume that
jobs arrive according to a general arrival process, and the job sizes
are independent of each other and identically distributed.

Note that letting $K=\infty$ makes the system a standard
\textit{processor sharing} (PS) queue, which has been the focal point of
significant recent research activity. The PS discipline can be viewed
as an idealization of time-sharing protocols in computer systems, as
described in \cite{Kleinrock1976} and \cite{RitchieThompson74}. The
advantage is that a big job will not block the whole system as in a
first-come-first-serve (FCFS) queue. However, allowing too many jobs
to time-share at once can lead to significant overhead due to
switching, and hence reduce overall performance. This point has
already been observed in early studies of operating systems
\cite{Blake82,Denningetal76}, as well as in more recent Web server
design papers \cite{Elniketyetal04,KamraMisraNahum04} and database
implementation papers \cite{HeissWagner91,ICDE06}. So in the modeling
of many computer and communication systems, a sharing limit is
normally imposed, which results in the LPS model.

Despite the numerous applications, there are only a few studies on the
LPS queue. Avi-Itzhak and Halfin \cite{Avi-ItzhakHalfin1988} propose
an approximation for the mean response time assuming Poisson arrivals.
A computational analysis based on matrix geometric methods is
performed in Zhang and Lipsky \cite{ZhaLi06a,ZhaLi07b}. Some
stochastic ordering results are derived in Nuyens and van de Weij
\cite{NuyensWeij07}. Recently, Zhang, Dai and Zwart \cite{ZDZ2009}
have developed a fluid approximation for the LPS queue using the
framework of measure-valued processes. As a continuation of
\cite{ZDZ2009}, the present study investigates a diffusion
approximation for the LPS queue in the heavy traffic regime.

In our model, the system consists of a server for serving jobs and a
buffer for holding the waiting jobs. We model the LPS queue by means
of a measure-valued process $(\buf(\cdot),\ser(\cdot))$. Each
component of the process takes values in the space of finite,
nonnegative Borel measures on $\R_+=[0,\infty)$. For each $t\ge0$,
$\ser(t)$ puts unit mass at the residual job size of each job in the
server at time $t\ge0$, and $\buf(t)$ puts unit mass at the job size
of each job in the buffer at time $t\ge0$. The main insight of our
approach is to design the stochastic dynamic equations
(\ref{eq:stoc-dym-eqn-B}) and (\ref{eq:stoc-dym-eqn-S}) using the
measure valued process $(\buf(\cdot),\ser(\cdot))$, which can describe
the evolution of the system. Our asymptotic regime is when the
sharing limit $K$ is large and the queue is critically
loaded. Following the standard practice in the literature, we consider
a sequence of queues indexed by $r\in\R_+$. We assume that
$\lim_{r\to\infty}K^r/r=K>0$ and the traffic intensity goes to the
critical value $1$ as in (\ref{eq:cond-HT}), where $K^r$ is the
sharing limit in the $r$th queue. (Superscript $r$ indicates a
quantity that is associated with the $r$th queue.)

We are interested in the limit of the diffusion scaled process
\[
\biggl(\frac{1}{r}\buf^r(r^2\cdot),\frac{1}{r}\ser^r(r^2\cdot
)\biggr)
\]
as $r$ goes to infinity. As shown in Williams \cite{Williams1998}, a
key step to obtain a diffusion limit in heavy traffic is to establish
a \textit{state-space collapse} (SSC) result. In our setting, the SSC
means that the diffusion-scaled measure-valued process, which is an
infinite-dimensional object, is close to a deterministic function of
the diffusion-scaled, one-dimensional workload process. (See
Definition \ref{defn:lifting-map} for the lifting map to define the
function.) The workload process is invariant under any nonidling
service policy, and its diffusion limit is a one-dimensional reflected
Brownian motion (RBM). The main result of this paper (Theorem
\ref{thm:diffusion}) is that our measure-valued diffusion limit is a
deterministic function of the one-dimensional RBM. As a corollary, the
diffusion-scaled system size process converges in distribution to a
piecewise linear RBM.

Most of this paper is devoted to proving the SSC result for each fixed
time $T>0$ (Theorem \ref{thm:ssc}). Our proof strategy is analogous
to the modular approach proposed in Bramson \cite{Bramson1998} and
Williams \cite{Williams1998}. For our sequence of systems, we define
a critically loaded measure-valued fluid model. The fluid model in
this paper is the same LPS fluid model developed in \cite{ZDZ2009},
specialized for the critically loaded case. We show that our fluid
model exhibits an SSC: each fluid model solution converges to an
equilibrium state in some uniform sense, and each equilibrium state
has an SSC.

We adopt Bramson's framework in \cite{Bramson1998} to translate the
fluid model SSC result into the diffusion-scaled SSC result. The
diffusion scaled process on the interval $[0,T]$ corresponds to
unscaled process on the interval $[0,r^2T]$. Fix a constant $L>1$,
the interval $[0, r^2T]$ is covered by the $\lfloor{rT}\rfloor+1$ overlapping
intervals
\[
[rm,rm+rL],\qquad m=0,1,\ldots,\lfloor{rT}\rfloor.
\]
On each of these
intervals, the diffusion scaled process can be viewed as a shifted,
fluid-scaled process defined by
\[
\biggl(\frac{1}{r}\buf^r(rm+rt), \frac{1}{r}\ser^r(rm+rt)\biggr),\qquad
0\le t \le L.
\]
To carry out the translation, we need to show that (a) each limit from
the family of shifted, fluid-scaled processes is a solution to the
fluid model (such a limit is called a fluid limit in this paper, also
known as a ``cluster point'' in \cite{Bramson1998}); (b) the set of
fluid limits is ``rich'': with large probability, each shifted,
fluid-scaled measure-valued process is close to some fluid limit. A
major step to proving (a) and (b) is to show, with large probability,
the precompactness of the shifted fluid scaled processes (see
Theorem \ref{thm:precompactness} for details). Since $m$ ranges from
$0$ to $\lfloor{rT}\rfloor$, this involves a substantial refinement
of the
arguments in \cite{ZDZ2009}, where the case $m=0$ is treated.

Establishing SSC for the fluid model requires a study of equilibrium
states for the LPS fluid model developed in \cite{ZDZ2009}. In
Section \ref{sec:fluid}, we characterize the set of equilibrium
states, and show that each fluid model solution
with initial condition belonging to a compact set
converges uniformly to its equilibrium state. The
counterpart of this study
for standard PS queues has been carried out in
\cite{PuhaWilliams2004}. Standard PS queues are relatively tractable
since their fluid models can be related (by means of a time-change) to
a renewal equation. This is not the case for LPS systems (as
explained in \cite{ZDZ2009}), so a different approach is necessary.
The idea behind the proof of the uniform convergence is to carefully
track the total mass of the fluid model, and determine whether it is
eventually bigger, smaller or equal to the sharing limit~$K$. Several
insights (explained in Sections \ref{subsec:conv-invar-manif} and
\ref{subsec:unif-conv-invar-manif}) lead us to apply a uniform version
of the renewal theorem which is new to the best of our knowledge. This
version is given in Appendix \ref{lem:key-renewal-thm-unif}.

The framework of measure-valued process has been successfully applied
to study models where multiple jobs are processed at the same time.
The main idea is to use a sufficiently detailed state descriptor to
adequately describe the system. A~sequence of papers, Gromoll, Puha
and Williams \cite{GPW2002}, Puha and Williams \cite{PuhaWilliams2004}
and Gromoll \cite{Gromoll2004}, has successfully established the fluid
and diffusion approximations for PS queues using measure-valued
processes. More recently, the framework of measure-valued process has
been further developed by Gromoll and Kruk \cite{GromollKurk2007} and
Gromoll, Robert and Zwart \cite{GRZ2008} in the study of queues with
deadlines/impatience. Doytchinov, Lehoczky and Shreve \cite{DLS2001}
applied a similar framework to study the earliest deadline first
discipline. This framework is also applied by Kaspi and Ramanan
\cite{KaspiRamanan} on many-server queues. The results in the present
paper can be seen as an extension of the results in the papers
\cite{Gromoll2004,PuhaWilliams2004}, which carry out a similar
program for the standard PS queue.

This paper is organized as follows. A model description and an
overview of the main results is given in
Section \ref{sec:models-main-result}. Section \ref{sec:fluid}
investigates the convergence of each fluid model solution to its
equilibrium state. Precompactness of the family of shifted fluid
scaled processes is established in Section \ref{sec:tightness}.
Section \ref{sec:diff-appr} uses the precompactness to show the
``richness'' of the fluid limits, and then concludes with a proof of
state-space collapse. Several additional useful results, such as a
uniform version of the renewal theorem, and a useful bound for the
Prohorov metric are developed in
Appendices \ref{sec:key-renewal-theorem} and \ref{sec:conv-proh-metr}.

\subsection{Notation}\label{subsec:notatioin}
The following notation will be used throughout. Let $\N$, $\Z$ and
$\R$ denote the set of natural numbers, integers and real numbers,
respectively. Let $\R_+=[0,\infty)$. For $a,b\in\R$, write $a^+$ for
the positive part of $a$, $\lfloor{a}\rfloor$ for the integer part,
$\lceil{a}\rceil$ for
$\lfloor{a}\rfloor+1$, $a\vee b$ for the maximum  and $a\wedge b$ for
the minimum.

Let $\M$ denote the set of all nonnegative finite Borel measures on
$[0,\infty)$. To simplify the notation, let us take the convention
that for any Borel set $A\subseteq\R$, $\nu(A\cap(-\infty,0))=0$ for
any $\nu\in\M$. For $\nu_1,\nu_2\in\M$, the Prohorov metric is
defined to be
\begin{eqnarray*}
\pov[\nu_1,\nu_2]&=&\inf\{\varepsilon>0\dvtx
\nu_1(A)\le\nu_2(A^\varepsilon)+\varepsilon\mbox{ and }\\
&&\hspace*{16.7pt}\nu_2(A)\le\nu_1(A^\varepsilon)+\varepsilon\mbox{ for all Borel set }
A\subseteq\R\},
\end{eqnarray*}
where $A^\varepsilon=\{b\in\R\dvtx{\inf_{a\in A}}|a-b|<\varepsilon\}$.
This same
metric was defined and used in Gromoll and Kruk
\cite{GromollKurk2007}; they showed that the space $\M$ is complete
and separable under the metric. For any Borel measurable function
$g\dvtx\R_+\to\R$, the integration of this function with respect to the
measure $\nu\in\M$, $\int_{\R_+}g(x)\nu(dx)$, is denoted by
$\langle{g},{\nu}\rangle$.

Let $\M\times\M$ denote the Cartesian product. There are a number of
ways to define the metric on the product space. For convenience, we
define the metric to be the maximum of the Prohorov metric between
each component. With a little abuse of notation, we still use $\pov$
to denote this metric.

Let $(\mathbf E, \gd)$ be a general metric space. We consider the
space $\D$ of all right-continuous $\mathbf E$-valued functions with
finite left limits defined either on a finite interval $[0,T]$ or the
infinite interval $[0,\infty)$. We refer to the space as
$\D([0,T],\mathbf E)$ or $\D([0,\infty),\mathbf E)$ depending on the
function domain. The space $\D$ is also known as the space of c\`adl\`ag
functions. For $g(\cdot),g'(\cdot)\in\D([0,T],\mathbf E)$, the
uniform metric is defined as
%
%
\begin{equation}\label{eq:sup-norm}
\upsilon_T[g,g']=\sup_{0\le t\le T}\mathbf\gd[g(t),g'(t)].
\end{equation}
However, a more useful metric we will use is the following Skorohod
$J_1$ metric:
%
%
\begin{equation}\label{eq:skorohod-L-def}
\sko_T[g,g']=\inf_{f\in\Lambda_T}(\|f\|^\circ_T\vee\upsilon
_T[g,g'\circ f]),
\end{equation}
where $g\circ f(\cdot)=g(f(\cdot))$ for $t\ge0$ and $\Lambda_T$ is
the set of strictly increasing and continuous mapping of $[0,T]$ onto
itself and
\[
\|f\|^\circ_{T}=\sup_{0\le s<t\le T}\biggl|\log\frac
{f(t)-f(s)}{t-s}\biggr|.
\]
If $g(\cdot)$ and $g'(\cdot)$ are in the space $\D([0,\infty
),\mathbf
E)$, the Skorohod $J_1$ metric is defined as
%
%
\begin{equation}\label{eq:skorohod-def}
\sko[g,g']=\int_0^\infty e^{-T}(\sko_T[g,g']\wedge1)\,dT.
\end{equation}
By ``convergence in the space $\D$,'' we mean the convergence under
the Skorohod $J_1$ topology, which is induced by the Skorohod $J_1$
metric \cite{EthierKurtz1986}.

We use ``$\to$'' to denote the convergence in the metric space
$(\mathbf E,\mathbf\gd)$, and ``$\dto$'' to denote the convergence in
distribution of random variables taking values in the metric space
$(\mathbf E,\mathbf\gd)$.

\section{Models and main results}\label{sec:models-main-result}
In this section, we first introduce the mathematical model. We then
present the main results of this paper. Following this, is an outline
of our proof.

\subsection{The limited processor sharing queue}
\label{subsec:limit-proc-shar-que}

We consider a $G/GI/1$ queue operated under the limited processor
sharing policy, with the sharing limit equal to~$K$. We use $Q(t)$,
$Z(t)$ and $X(t)$ to denote the number of jobs in the buffer, the
number of jobs in service, and the total number of jobs in the system
at time $t$, respectively. Thus,
%
%
\begin{equation}\label{eq:X(t)}
X(t)=Q(t)+Z(t) \qquad\mbox{for }t\ge0.
\end{equation}
We adopt the convention that $Q(\cdot),Z(\cdot)$ and $X(\cdot)$ are
right continuous. The system is allowed to be nonempty initially,
that is, $X(0)>0$. We index jobs by $i=-X(0)+1, -X(0)+2, \ldots, 0,
1,\ldots.$ The first $X(0)$ jobs are initially in the system, with
jobs $i=-X(0)+1, \ldots, -Q(0)$ in service and jobs $i=-Q(0)+1, \ldots,
0$ waiting in the buffer. Jobs arrived after time $0$ are indexed by
$i=1,2,\ldots,$ according to the order of arrival. When a batch arrival
occurs, an arbitrary rule is used to break the tie for the arrivals in
the batch.
The service policy in this model is FCFS. Let $E(t)$ denote the
number of jobs that arrive at the buffer during time interval $(0,t]$,
for all $t\ge0$. Our arrival process $\{E(t), t\ge0\}$ is assumed to
be general, as long as it satisfies a functional central limit theorem
[see (\ref{eq:cond-A})]. According to the policy, a job may have to
wait for a certain amount of time after arrival to get service. Let
$w_i$ denote the waiting time, and $U_i$ denote the arrival time of
the $i$th job for all $i>-X(0)$. By convention, $U_i=0$ for $i<0$,
and $w_i=0$ for $i\le-Q(0)$. Let
\[
\tau_i=U_i+w_i,\qquad i> -X(0).
\]
The quantity $\tau_i$ can be viewed as the time that the $i$th job
starts service. We use $v_i$ to denote the job size of the $i$th job
for all $i>-Q(0)$. We assume that $\{v_i\}_{i=-\infty}^\infty$ is a
sequence of i.i.d. random variables with distribution
$F(\cdot)$. Denote $\nu$ the probability measure associated with the
distribution function $F(\cdot)$. For jobs with index $-X(0)<i\le
-Q(0)$, that is, the first $Z(0)$ jobs that are initially in service, we
use $\tilde v_i$ to denote the remaining job size of the job. The
sequence $\{\tilde v_i\}_{i=-\infty}^0$ is allowed to be general. We
call $\{E(\cdot),\{v_i\}_{i=1}^\infty\}$ the stochastic primitives of
the system, and $\{Z(0),Q(0),\{v_i\}_{i=-\infty}^0,\{\tilde
v_i\}_{i=-\infty}^0\}$ the initial conditions of the system.

Now we introduce a measure-valued state descriptor
$(\buf(\cdot),\ser(\cdot)) \in\M\times\M$, which describes the
evolution of the system with given initial conditions and stochastic
primitives. For any Borel set $A\subset[0,\infty)$, $\buf(t)(A)$
denotes the total number of jobs in the buffer whose job size belongs
to $A$; and for any Borel set $A\subset(0,\infty)$, $\ser(t)(A)$
denotes the total number of jobs in service whose residual job size
belongs to $A$. Since no job can be in service with residual job size
$0$, $\ser(t)(\{0\})=0$ for all $t\ge0$. It is clear that we have the
following relationship:
\[
Q(t)=\langle{1},{\buf(t)}\rangle,\qquad Z(t)=\langle{1},{\ser
(t)}\rangle.
\]
Define the \textit{cumulative service amount} up to time $t$ by
%
%
\begin{equation}\label{eq:cumulative}
S(t)=\int_0^t\psi(Z(\tau))\,d\tau,
\end{equation}
where $\psi(x)=1/x$ if $x>0$ and $\psi(x)=0$ if $x=0$. A job will
have received a cumulative amount of processing time
\[
S(s,t)=\int_s^t\psi(Z(\tau))\,d\tau
\]
during time interval $[s,t]$ if it is in service in this time period.
Let
%
%
\begin{equation}\label{eq:B(t)}
B(t)=E(t)-Q(t).
\end{equation}
Note that at time $t\ge0$, $B(t)$ is the index of the last job which
has entered into service by time $t$. Thus,
%
%
\begin{equation}
B(s,t)=B(t)-B(s)
\end{equation}
represents the number of jobs which have left the buffer and entered
the server during time interval $(s,t]$. Using the notation
introduced in this section, the state descriptor can be written as
%
%
\begin{eqnarray}
\label{eq:stoc-dym-eqn-B}
\buf(t)(A)&=&\sum_{i=B(t)+1}^{E(t)}\delta_{v_i}(A),\qquad A\subset
[0,\infty),\\
\label{eq:stoc-dym-eqn-S}
\ser(t)(A)&=&\sum_{i=-X(0)+1}^{-Q(0)}\delta_{\tilde
v_i}\bigl(A+S(t)\bigr)\nonumber\\[-8pt]\\[-8pt]
&&{}+\sum_{i=- Q(0)+1}^{B(t)}\delta_{v_i}\bigl(A+S(\tau_i,t)\bigr),\qquad A\subset
(0,\infty),\nonumber
\end{eqnarray}
with $\ser(t)(\{0\})=0$, where $\delta_a(A)$ denotes the Dirac measure
of point $a$ on $\R$ and $A+y=\{a+y\dvtx a\in A\}$. Due to the LPS policy,
the sharing limit $K$ must be enforced at any time $t$,
%
%
\begin{eqnarray}
\label{eq:policy-B}Q(t)&=&\bigl(X(t)-K\bigr)^+,\\
\label{eq:policy-S}Z(t)&=&\bigl(X(t)\wedge K\bigr).
\end{eqnarray}
We call (\ref{eq:stoc-dym-eqn-B}) and (\ref{eq:stoc-dym-eqn-S}) the
\textit{stochastic dynamic equations} and (\ref{eq:policy-B}) and
(\ref{eq:policy-S}) the policy constraints.

For $t\ge0$, the workload of the system $W(t)$ is defined to be the
amount of time that the server remains busy \textit{if} no more arrivals
are allowed into the system at time~$t$. Using the state descriptor
$(\buf,\ser)$, we can recover the workload $W(t)$ at time $t\ge0$ by
%
%
\begin{equation}
W(t)=\langle{\chi},{\buf(t)+\ser(t)}\rangle,
\end{equation}
where $\chi$ denotes the identity function on $\R$.

\subsection{Main results}\label{subsec:main-results}

Consider a sequence of limited processor sharing queues indexed by
$r$, where $r$ increases to $\infty$ through a sequence in
$(0,\infty)$. Each queue is defined in the same way as in
Section \ref{subsec:limit-proc-shar-que}. To distinguish models with
different indices, quantities of the $r$th model are accompanied by
superscript $r$. Each model may be defined on a different probability
space $(\Omega^r,\mathcal F^r,\mathbb P^r)$. Our results concern the
asymptotic behavior of the descriptor under the \textit{diffusion}
scaling, which is defined by
%
%
\begin{equation}
\dr\buf(t)=\frac{1}{r} \buf^r(r^2t),\qquad
\dr\ser(t)=\frac{1}{r} \ser^r(r^2t),
\end{equation}
for all $t\ge0$. We are also interested in other diffusion scaled
quantities like the workload and queue length processes. Note that
$Q^r(\cdot)$, $Z^r(\cdot)$ and $W^r(\cdot)$ are actually functions of
$(\buf^r(\cdot),\ser^r(\cdot))$, so the scaling for these quantities
is defined as the functions of the corresponding scaling for
$(\buf^r(\cdot),\ser^r(\cdot))$, that is,
%
%
\begin{eqnarray}
\dr{Q}(t)&=&\langle{1},{\dr\buf(t)}\rangle=\frac{1}{r}Q^r(r^2t),\\
\dr{Z}(t)&=&\langle{1},{\dr\ser(t)}\rangle=\frac{1}{r}Z^r(r^2t),\\
\dr{W}(t)&=&\langle{\chi},{\dr\buf(t)+\dr\ser(t)}\rangle=\frac
{1}{r}W^r(r^2t),
\end{eqnarray}
for all $t\ge0$.

To establish results on the convergence of the above sequence of
stochastic processes, we need the following conditions, which are
quite general and standard. We assume that the arrival processes
satisfy
%
%
\begin{equation}\label{eq:cond-A}
\frac{E^r(r^2\cdot)-\lambda^r r^2 \cdot}{r}\dto E^*(\cdot)\qquad
\mbox{as }r\to\infty,
\end{equation}
for some sequence $\{\lambda^r\}$ that satisfies
%
%
\begin{equation}\label{eq:cond-lambda}
\lim_{r\to\infty}\lambda^r=\lambda>0,
\end{equation}
and $E^*(\cdot)$ is a Brownian motion with drift $0$ and variance
$\lambda c_a^2$. And the probability measure $\nu^r$ of job sizes in
the $r$th system satisfies that as $r\to\infty$
%
%
\begin{eqnarray}
\label{eq:cond-S}\pov[\nu^r,\nu]&\to&0,\\
\label{eq:cond-S-u}\lim_{N\to\infty}\sup_{r}
\int_{[N, \infty)} x^{4+2p} \nu^r(dx)&\to&0
\qquad\mbox{for some }p>0,
\end{eqnarray}
where the probability measure $\nu$ satisfies
%
%
\begin{equation}
\label{eq:cond-nu-noatom}
\nu\mbox{ has no atoms}.
\end{equation}
%
Assumption (\ref{eq:cond-S-u}) is stronger than the ``two-plus-epsilon
moment'' assumption needed for a functional central limit theorem. The
stronger assumption is used in a separate part of our analysis to
estimate moments of the shifted fluid scaled state descriptors (see
Lemma \ref{lem:glivenko-cantelli}, and its application in
Lemma \ref{lem:p-moment-bound}). The extra moment assumption also
appears in \cite{Gromoll2004,GromollKurk2007}.
Since the space is scaled by $r$ in the diffusion scaling, the sharing
limit should be scaled accordingly:
%
%
\begin{equation}\label{eq:cond-K}
\lim_{r\to\infty}K^r/r\to K>0.
\end{equation}
Also, the following initial condition will be assumed:
%
%
\begin{eqnarray}
\label{eq:cond-initial}
(\dr\buf(0),\dr\ser(0))&\dto&(\xi^*,\mu^*),\\
\label{eq:cond-initial-u}
\langle{\chi^{1+p}},{\dr\buf(0)+\dr\ser(0)}\rangle&\dto&\langle
{\chi^{1+p}},{\xi^*+\mu^*}\rangle,
\end{eqnarray}
as $r\to\infty$, where $p$ is the same as in (\ref{eq:cond-S-u}),
$(\xi^*,\mu^*)\in\M\times\M$ and
%
%
\begin{equation}\label{eq:cond-initial-noatom}
\mu^* \mbox{ has no atoms}.
\end{equation}
%
Define the traffic intensity of the $r$th stochastic system by
$\rho^r=\lambda^r\langle{\chi},{\nu^r}\rangle$. We need the
following heavy
traffic condition:
%
%
\begin{equation}\label{eq:cond-HT}
\lim_{r\to\infty}r(1-\rho^r)=\theta>0.
\end{equation}

Let $\beta=\langle{\chi},{\nu}\rangle$ be the mean and
$c^2_s=\frac{\langle{\chi^2},{\nu}\rangle-\beta^2}{\beta^2}$ be
the squared
coefficient of variation (SCV) of the job size distribution $\nu$.
The following proposition is a well-known heavy traffic approximation
for the workload process of a single queue operated under a nonidling
policy. Readers are referred to \cite{Gromoll2004} for a proof.
\begin{prop}\label{prop:workload}
Assume (\ref{eq:cond-A})--(\ref{eq:cond-S-u}),
(\ref{eq:cond-initial}), (\ref{eq:cond-initial-u}) and
(\ref{eq:cond-HT}). The sequence of diffusion scaled workload
process
\[
\dr W(\cdot)\dto W^*(\cdot) \qquad\mbox{as }r\to\infty,
\]
where $W^*(\cdot)$ is a reflected Brownian motion with drift
$-\theta$, variance $\beta(c_a^2+c_s^2)$ and initial value
$w^*=\langle{\chi},{\xi^*+\mu^*}\rangle$.
\end{prop}

Since the LPS is also a nonidling service policy, the above result on
the workload process is still true for our model. However, it remains
an open question about the job size process $X(\cdot)$ and many other
performance processes as introduced in Section
\ref{subsec:limit-proc-shar-que}. Our main result establishes the
diffusion limit for the measure-valued processes
(Theorem \ref{thm:diffusion}), from which the diffusion limit of queue
length process follows directly (Corollary \ref{cor:piecewise-RBM}).

Denote
\[
\beta_e=\langle{\chi},{\nu_e}\rangle,
\]
where $\nu_e$ is the \textit{equilibrium} measure of $\nu$, that is,
$\nu_e([0,x])=\frac{1}{\beta}\int_0^x\nu((y,\infty))\,dy$ for all
$x\ge
0$. We have the following definition.
\begin{defn}\label{defn:lifting-map}
Let $\Delta_{K,\nu}\dvtx\R_+\to\M\times\M$ be the lifting map associated
with the probability measure $\nu$ and constant $K$ given by
\[
\Delta_{K,\nu} w=\biggl(\frac{(w-K\beta_e)^+}{\beta}\nu,
\frac{w\wedge K\beta_e}{\beta_e}\nu_e\biggr)\qquad
\mbox{for }w\in\R_+.
\]
\end{defn}

Note that $\Delta_{K,\nu}$ maps the workload, which is in $\R_+$, to a
measure-valued state, which is in $\M\times\M$. The intuition is that
the remaining job sizes of those in service have the probability
measure $\nu_e$, which is the equilibrium distribution of the job size
distribution $\nu$. The total workload embodied in jobs that are in
service is $w\wedge K\beta_e$ because at most $K$ jobs are allowed in
service. Dividing it by $\beta_e$ gives the number of jobs in
service. The remaining workload $(w-K\beta_e)^+=w-(w\wedge K\beta_e)$
resides in the buffer, where each job size follows the probability
measure $\nu$. Dividing this amount by $\beta$, the mean of $\nu$,
gives the number of jobs in buffer.

Our main result requires that the limit $(\xi^*,\mu^*)$ in
(\ref{eq:cond-initial}) satisfies
%
%
\begin{equation}\label{eq:cond-init-ssc}
(\xi^*,\mu^*)=\Delta_{K,\nu}w^*.
\end{equation}
A similar condition on the initial state is also imposed for proving a
diffusion limit in queueing networks
\cite{Bramson1998,Williams1998}. When (\ref{eq:cond-init-ssc}) does
not hold, it is likely that some delayed versions of theorems are
still true as in Theorem 3 of Bramson \cite{Bramson1998}; however, we
will not pursue such generalization in this paper.
\begin{theorem}\label{thm:diffusion}
Assume (\ref{eq:cond-A})--(\ref{eq:cond-HT}). The sequence of
diffusion scaled state descriptors
\[
(\dr\buf(\cdot),\dr\ser(\cdot))\dto
\Delta_{K,\nu}W^*(\cdot) \qquad\mbox{as }r\to\infty,
\]
where $W^*(\cdot)$ is the
reflected Brownian motion in Proposition \ref{prop:workload}.
\end{theorem}
\begin{cor}[(Piecewise reflected Brownian motion)]\label{cor:piecewise-RBM}
Assume (\ref{eq:cond-A})--(\ref{eq:cond-HT}). The sequence of
diffusion scaled system size process
\[
\dr X(\cdot)=\langle{1},{\dr\buf(\cdot)+\dr\ser(\cdot)}\rangle
\]
converges in distribution as $r\to\infty$ to $X^*(\cdot)$, where
\[
X^*(t)=\frac{(W^*(t)-K\beta_e)^+}{\beta}
+\frac{W^*(t)\wedge K\beta_e}{\beta_e}\qquad
\mbox{for }t\ge0,
\]
and $W^*(\cdot)$ is the reflected Brownian motion as in
Proposition \ref{prop:workload}.
\end{cor}
\begin{pf}
Since $\hat X^r(\cdot)=\langle{1},{\hat\buf^r(\cdot)+\hat\ser
^r(\cdot)}\rangle$
and the mapping $\Phi\dvtx\M\times\M\to\R$ defined by
$\Phi(\nu_1,\nu_2)=\langle{1},{\nu_1+\nu_2}\rangle$ for any
$(\nu_1,\nu_1)\in\M\times\M$ is continuous, the result follows from
Theorem \ref{thm:diffusion} and the continuous mapping theorem.
\end{pf}
\begin{rem}
In other words, $X^*(\cdot)$ is a reflected Brownian motion with
drift $\frac{-\theta}{\beta}$ and variance
$\frac{c_a^2+c_s^2}{\beta^2}$ when it is above $K\beta_e$ and with
drift $\frac{-\theta}{\beta_e}$ and variance
$\frac{c_a^2+c_s^2}{\beta_e^2}$ when it is below $K\beta_e$.
\end{rem}

\subsection{Outline of proof}

The major step to prove our main result is to establish the following
\textit{state-space collapse} result.
\begin{theorem}\label{thm:ssc}
Assume (\ref{eq:cond-A})--(\ref{eq:cond-HT}). Fix $T>0$. If
(\ref{eq:cond-init-ssc}) holds, then
\[
\sup_{t\in[0,T]}\pov[(\dr\buf(t),\dr\ser(t)),\Delta_{K,\nu}\dr W(t)]
\dto0 \qquad\mbox{as }r\to\infty.
\]
\end{theorem}

The state-space collapse result is appealing, since it rigorously
shows that all performance processes can be described as a simple,
deterministic function of the workload process. We now use
Theorem \ref{thm:ssc} to prove the main result.
\begin{pf*}{Proof of Theorem \ref{thm:diffusion}}
We have the convergence of the workload pro\-cesses $\dr W(t)$ in
Proposition \ref{prop:workload}. Since the mapping
$\Delta_{K,\nu}\dvtx\R_+\to\M\times\M$ is continuous, by the continuous
mapping theorem
\[
\Delta_{K,\nu}\dr W(\cdot)\dto\Delta_{K,\nu}W^*(\cdot)\qquad
\mbox{as }r\to\infty.
\]
The result of the theorem follows immediately from the state-space
collapse result in Theorem \ref{thm:ssc} and the ``convergence
together lemma'' (Theorem 4.1 in \cite{Billingsley1999}).
\end{pf*}

The proof of the state-space collapse is given in
Section \ref{sec:diff-appr}, but it requires ample preparation. Our
proof is analogous to the framework developed in
Bramson \cite{Bramson1998} for proving state-space collapse in
multi-class networks with head-of-the-line service disciplines. The
framework was later adopted for the PS queue in
Gromoll~\cite{Gromoll2004}, which the current paper closely follows.
In Section \ref{sec:fluid}, we establish several fundamental
properties for the equilibrium behavior of the LPS fluid model
introduced in \cite{ZDZ2009}. Section \ref{sec:tightness} establishes
precompactness of a family of shifted fluid scaled processes, which
will be defined in that section. Briefly speaking, the proof of the
state-space collapse is built on the ``richness'' of the set of fluid
limits, which are obtained from the shifted fluid scaled processes.
Each fluid limit is shown to be a fluid model solution in
Section \ref{subsec:fluid-limits}.

\section{Convergence to equilibrium states for fluid model}
\label{sec:fluid}

We propose a fluid model, denoted by $(K,\lambda, \nu)$, to assist the
study of the underlying stochastic processes for the LPS queue. The
parameters $K$, $\lambda$ and $\nu$ are the limiting sharing level
defined in (\ref{eq:cond-K}), the limiting arrival rate defined in
(\ref{eq:cond-lambda}), and the limiting job size distribution defined
(\ref{eq:cond-S}), respectively. According to condition
(\ref{eq:cond-lambda}), (\ref{eq:cond-S}), (\ref{eq:cond-S-u}) and
(\ref{eq:cond-HT}), we have that the traffic intensity of the fluid
model
%
%
\begin{equation}\label{eq:critical-load-rho}
\rho=\lambda\beta=1.
\end{equation}
Although $\nu$ is required to satisfy (\ref{eq:cond-nu-noatom}) in
Theorem \ref{thm:diffusion}, in this section only, we allow the job
size distribution to have atoms.
The fluid analogue of the LPS queue was first proposed and studied in
\cite{ZDZ2009}, where general properties of the fluid model were
studied for all traffic intensities $\rho\in[0,\infty)$. For the
purpose of this paper, we now recall the definition and some general
properties for critically loaded fluid model.

Given a measure-valued process
$(\bar\buf(\cdot),\bar\ser(\cdot))\in\D([0,\infty),\M\times
\M)$, for
$t\ge0$, let
%
%
\begin{eqnarray}\label{eq:fluid-def-Q}
\bar Q(t)&=&\langle{1},{\bar\buf(t)}\rangle,\\
\bar Z(t)&=&\langle{1},{\bar\ser(t)}\rangle,\\
\bar X(t)&=&\bar Q(t)+\bar Z(t),\\
\bar B(t)&=&\lambda t-\bar Q(t).
\end{eqnarray}
These quantities are the fluid analogues of $Q(t), Z(t), B(t)$ and
$X(t)$ in the stochastic model. Define the \textit{fluid cumulative
service amount} up to time $t$ by
%
%
\begin{equation}\label{eq:fluid-accumulate}
\bar S(t)=\int_0^t\phi_\rho(\bar Z(\tau))\,d\tau,
\end{equation}
where
%
%
\begin{equation}\label{eq:def-psi_rho}
\phi_\rho(x)=\cases{
1/x, &\quad $x>0$,\cr
\infty, &\quad $x=0$,}
\end{equation}
when $\rho=1$. It worth noting that the function
$\phi_\rho$ take a slightly different form when $\rho\ne1$, which is
not the case in this paper. Interested readers are referred to
\cite{ZDZ2009} for detailed discussion.
And for $0\le s\le t$, denote
%
%
\begin{equation}\label{eq:fluid-accumulate-interval}
\bar S(s,t)=\int_s^t\phi_\rho(\bar Z(\tau))\,d\tau.
\end{equation}

An element $(\xi,\mu)\in\M\times\M$ is called a \textit{valid initial
condition} if
%
%
\begin{eqnarray}\label{eq:valid-init-1}
{\xi}&=&(\langle{1},{\xi}\rangle+\langle{1},{\mu}\rangle-K)^+\nu
,\\
\label{eq:valid-init-2}
\langle{1},{\mu}\rangle&=&(\langle{1},{\xi}\rangle+\langle
{1},{\mu}\rangle)\wedge K.
\end{eqnarray}
Roughly speaking, validity of an initial state means that the initial
state is consistent with the LPS policy; initial waiting jobs have the
same job size distribution as arriving jobs. Denote
%
%
\begin{equation}\label{eq:valid-init}
\init=\{(\xi,\mu)\in\M\times\M\dvtx(\xi,\mu)
\mbox{ satisfies (\ref{eq:valid-init-1})
and (\ref{eq:valid-init-2})}\}
\end{equation}
the set of all valid initial conditions.

We now introduce the following \textit{fluid dynamic equations}, which
are analogous to (\ref{eq:stoc-dym-eqn-B}) and
(\ref{eq:stoc-dym-eqn-S}). For all $A_y=(y,\infty)$, $y\ge0$,
%
%
\begin{eqnarray}
\label{eq:fluid-dym-eqn-B}
\bar\buf(t)(A_y)&=&\xi(A_y)+\bigl(\bar Q(t)-\bar Q(0)\bigr)\nu
(A_y),\\
\label{eq:fluid-dym-eqn-S}
\bar\ser(t)(A_y)&=&\mu\bigl(A_y+\bar S(t)\bigr)
+\int_{0}^{t}\nu\bigl(A_y+\bar S(s,t)\bigr)\,d\bar B(s),
\end{eqnarray}
where $\bar Q(\cdot)$, $\bar Z(\cdot)$, $\bar X(\cdot)$, $\bar
B(\cdot)$ and $\bar S(\cdot)$ are defined in
(\ref{eq:fluid-def-Q})--(\ref{eq:fluid-accumulate-interval}). They
are subject to the following constraints:
%
%
\begin{eqnarray}
\label{constr:non-decr}
&&\bar B(\cdot) \mbox{ is nondecreasing,}\\
\label{constr:policy-B}
&&\bar Q(t)=\bigl(\bar X(t)-K\bigr)^+,\\
\label{constr:policy-S}
&&\bar Z(t)=\bigl(\bar X(t)\wedge K\bigr).
\end{eqnarray}
Because
$(\bar\buf(\cdot),\bar\ser(\cdot))\in\D([0,\infty),\M\times
\M)$, $\bar
Q(\cdot)$, $\bar X(\cdot)$, $\bar Z(\cdot)$ and $\bar B(\cdot)$ are
right continuous on $[0,\infty)$ and have left limits in
$(0,\infty)$. Here and later, the integral $\int_0^tg(s)\,d\bar B(s)$ is
interpreted as the Lebesgue--Stieltjes integral on the interval
$[0,t]$, where by convention we set $B(0^-)=B(0)$.
The above dynamic equations and constraints define the fluid model
$(K,\lambda,\nu)$.
\begin{defn}\label{def:fluid}
$(\bar\buf(\cdot),\bar\ser(\cdot))\in\D([0,\infty),\M\times
\M)$ is a
solution to the fluid model $(K,\lambda,\nu)$ with a valid initial
condition $(\xi,\mu)$ if it satisfies the fluid dynamic equations
(\ref{eq:fluid-dym-eqn-B}) and (\ref{eq:fluid-dym-eqn-S}), subject
to the constraints (\ref{constr:non-decr})--(\ref{constr:policy-S}).
\end{defn}


Note that (\ref{eq:cond-S-u}) and (\ref{eq:cond-S-u}) imply that
%
%
\begin{equation}\label{eq:cond-nu-1}
\beta<\infty,
\end{equation}
and (\ref{eq:cond-nu-noatom}) implies
%
%
\begin{eqnarray}
\label{eq:cond-nu-11}&&\nu(\{0\})=0,\\
\label{eq:cond-nu-nonlattice}&&\nu\mbox{ is nonlattice}.
\end{eqnarray}
It has been proved in \cite{ZDZ2009} (cf. Theorem 3.1) that under
the conditions (\ref{eq:cond-nu-1}) and (\ref{eq:cond-nu-11}), there
exists a unique fluid model solution
$(\bar\buf(\cdot),\bar\ser(\cdot))$ for any valid initial condition
$(\xi,\mu)$. Moreover, by Proposition 3.1 in \cite{ZDZ2009}, the
fluid workload process which is defined as
$\bar W(t)=\langle{\chi},{\bar\buf(t)+\bar\ser(t)}\rangle$
satisfies the workload conservation property, that is,
\[
\bar W(t)=\bigl(\langle{\chi},{\xi+\mu}\rangle+(\rho-1)t\bigr)^+
\qquad\mbox{for all } t\ge0.
\]
Since we restrict to the critically loaded case, that is, $\rho=1$, we
have
%
%
\begin{equation}\label{eq:prop-workload}
\bar W(t)=\langle{\chi},{\xi+\mu}\rangle\qquad\mbox{for all
}t\ge0.
\end{equation}
The main objective of this section is to show the following long-term
behavior of the critically loaded fluid model, which helps to
establish the state-space collapse in Section \ref{sec:diff-appr}.
\begin{theorem}\label{thm:convergence-to-invariant-mani}
Assume (\ref{eq:critical-load-rho}) and
(\ref{eq:cond-nu-1})--(\ref{eq:cond-nu-nonlattice}). The unique
solution $(\bar\buf(\cdot)$, $\bar\ser(\cdot))$ to the fluid model
$(K,\lambda,\nu)$ with a valid initial state $(\xi,\mu)$ such that
$w=\langle{\chi},{\xi+\mu}\rangle<\infty$ satisfies
\[
(\bar\buf(t),\bar\ser(t))\to\Delta_{K,\nu} w \qquad\mbox{as }
t\to\infty.
\]
Moreover, for fixed constants $p, M>0$ the convergence is uniform
for all fluid solutions with initial conditions in the set
%
%
\begin{equation}\label{eq:comp-set-init}
\init_M^p=\{(\xi,\mu)\in\init\dvtx\langle{\chi},{\xi+\mu}\rangle<M,
\langle{\chi^{1+p}},{\xi+\mu}\rangle<M\}.
\end{equation}
\end{theorem}

Section \ref{subsec:invariantmanifold} characterizes the equilibrium
states for the fluid model. Section \ref{subsec:conv-invar-manif}
presents the proof of convergence (the first half of
Theorem \ref{thm:convergence-to-invariant-mani}), and
Section \ref{subsec:unif-conv-invar-manif} presents the proof of
uniform convergence (the second half of
Theorem \ref{thm:convergence-to-invariant-mani}).

\subsection{Equilibrium states}\label{subsec:invariantmanifold}
%
\begin{defn}\label{def:invariant-manifold}
An element $(\xi,\mu)\in\init$ is called an \textit{equilibrium state}
for the fluid model $(K,\lambda,\nu)$ if the solution to the fluid
model with initial condition $(\xi,\mu)$ satisfies
\[
(\bar\buf(t),\bar\ser(t))=(\xi,\mu) \qquad\mbox{for all }t\ge0.
\]
\end{defn}

The simple intuition here is that if the fluid model solution starts
with an invariant state, it will stay at this state forever. By the
restarting lemma, Lemma 4.2 in~\cite{ZDZ2009}, a fluid model solution
will remain in an equilibrium state once reach it. Our first result is
a characterization of an equilibrium state.
\begin{theorem}\label{thm:invariant-manifold}
An element $(\xi,\mu)\in\init$ is an equilibrium state if and only
if
%
%
\begin{equation}\label{eq:inv-mani}
(\xi,\mu)=\Delta_{K,\nu} w \qquad\mbox{for some }w\in[0,\infty).
\end{equation}
\end{theorem}
\begin{pf}
Suppose $(\xi,\mu)=\Delta_{K,\nu} w$ for some $w\in[0,\infty)$, we
need show that
\[
(\bar\buf(\cdot),\bar\ser(\cdot))\equiv\Delta_{k,\nu}w
=\biggl(\frac{(w-K\beta_e)^+}{\beta}\nu,\frac{w\wedge K\beta_e}{\beta
_e}\nu_e\biggr)
\]
is the fluid model solution. If $w=0$, then by weak stability
(Theorem 3.2 in \cite{ZDZ2009}), $\Delta_{k,\nu}0=(\mathbf
{0},\mathbf{0})$ is the
fluid model solution. So let us now focus on the case where $w> 0$.
The fluid amount of jobs in buffer size and in service are
\begin{eqnarray*}
\bar Q(t)&=&\langle{1},{\bar\buf(t)}\rangle=\frac{(w-K\beta
_e)^+}{\beta},\\
\bar Z(t)&=&\langle{1},{\bar\ser(t)}\rangle=\frac{w\wedge K\beta
_e}{\beta_e}.
\end{eqnarray*}
If $\bar Z(t)<K$, then $w<K\beta_e$ which implies that $\bar
Q(t)=0$; if $\bar Q(t)>0$, then $w>K\beta_e$ which implies that
$\bar Z(t)=K$. So condition (\ref{constr:policy-B}) and
(\ref{constr:policy-S}) in Definition \ref{def:fluid} are satisfied.
Since $\bar Q(t)$ and $\bar Z(t)$ remain to be a constant,
(\ref{constr:non-decr}) holds trivially. This also implies that the
fluid dynamic equation (\ref{eq:fluid-dym-eqn-B}) is satisfied. It
remains to verify the fluid dynamic equation
(\ref{eq:fluid-dym-eqn-S}). The fluid accumulative service amount
\[
\bar S(t)=\frac{t}{\bar Z(0)}=\frac{\beta_e}{w\wedge K\beta_e}t
\]
since $\bar Z(t)$ is a constant. The right-hand side of
(\ref{eq:fluid-dym-eqn-S}) becomes
\[
\frac{w\wedge K\beta_e}{\beta_e}
\nu_e\biggl(A_y+\frac{\beta_e}{w\wedge K\beta_e}t\biggr)
+\lambda\int_{0}^{t}\nu\biggl(A_y+\frac{\beta_e}{w\wedge K\beta_e}(t-s)\biggr)\,ds,
\]
which equals $\frac{w\wedge
K\beta_e}{\beta_e}\nu_e(A_y)=\bar\ser(t)(A_y)$, for all $y\ge0$.
So (\ref{eq:fluid-dym-eqn-S}) is verified. Thus,
$(\bar\buf(\cdot),\bar\ser(\cdot))$ is the fluid model solution.

Suppose that $(\xi,\mu)$ is an equilibrium state, we need to show
that $(\xi,\mu)$ takes the form (\ref{eq:inv-mani}). If
$(\xi,\mu)=(\mathbf{0},\mathbf{0})$, then trivially $(\xi,\mu
)=\Delta_{K,\nu} 0$.
Let us now assume that $(\xi,\mu)\ne(\mathbf{0},\mathbf{0})$. Since
$(\bar\buf(\cdot),\bar\ser(\cdot))\equiv(\xi,\mu)$ is the
fluid model
solution, the fluid dynamic equation (\ref{eq:fluid-dym-eqn-S}) must
be satisfied. That is,
\begin{eqnarray*}
\mu(A_y)&=&\mu\biggl(A_y+\frac{t}{\langle{1},{\mu}\rangle}\biggr)
+\lambda\int_{0}^{t}\nu\biggl(A_y+\frac{t-s}{\langle{1},{\mu}\rangle
}\biggr)\,ds,\\
&=&\mu\biggl(A_y+\frac{t}{\langle{1},{\mu}\rangle}\biggr)
+\frac{\langle{1},{\mu}\rangle}{\beta}\int_y^{y+{t}/{\langle
{1},{\mu}\rangle}}\nu
(A_{s'})\,ds'
\end{eqnarray*}
for all $y,t\ge0$. This yields
\[
\mu(A_y)-\mu\biggl(A_y+\frac{t}{\langle{1},{\mu}\rangle}\biggr)=
\langle{1},{\mu}\rangle\biggl(\nu_e(A_y)-\nu_e\biggl(A_y+\frac
{t}{\langle{1},{\mu}\rangle}\biggr)\biggr),
\]
which implies that $\mu=\langle{1},{\mu}\rangle\nu_e$ due to the
arbitrary of
$t$ and $y$. Since $(\xi,\mu)$ is a valid state,
$\xi=\langle{1},{\xi}\rangle\nu$. Let
\[
w=\langle{\chi},{\xi+\mu}\rangle=\langle{1},{\xi}\rangle\beta
+\langle{1},{\mu}\rangle\beta_e.
\]
Again by validity of state $(\xi,\mu)$,
$\langle{1},{\xi}\rangle=\frac{(w-K\beta_e)^+}{\beta}$ and
$\langle{1},{\mu}\rangle=\frac{w\wedge K\beta_e}{\beta_e}$. So we conclude
that $(\xi,\mu)=\Delta_{K,\nu} w$.
\end{pf}

\subsection{Convergence to equilibrium states}
\label{subsec:conv-invar-manif}

We now identify conditions under which the fluid model solution
starting at a valid initial state $(\xi,\mu)$ will converge to an
equilibrium state.

If the initial condition $(\xi, \mu)=(\mathbf{0},\mathbf{0})$, then
by weak stability
(Theorem 3.2 in \cite{ZDZ2009}), the fluid model solution will
always be zero. So $(\mathbf{0},\mathbf{0})$ is an equilibrium state.
From now on,
we focus on the case where the initial condition $(\xi,
\mu)\ne(\mathbf{0},\mathbf{0})$. By the fluid dynamic equation
(\ref{eq:fluid-dym-eqn-B}), $\bar Q(t)\beta=\langle{\chi},{\bar
\buf(t)}\rangle$.
It follows from the workload conservation property
(\ref{eq:prop-workload}) that $w=\langle{\chi},{\xi+\mu}\rangle
\equiv\bar
W(t)\ge\langle{\chi},{\bar\buf(t)}\rangle$ for all $t\ge0$. So
%
%
\begin{equation}\label{eq:upper-bound}
\bar Q(t)=\bigl(\bar X(t)-K\bigr)^+\le\frac{w}{\beta} \qquad\mbox{for all }
t\ge0.
\end{equation}
Since $\bar W(t)=0$ if and only if $\bar Z(t)=0$,
%
%
\begin{equation}
\bar Z(t)=\bigl(\bar X(t)\wedge K\bigr)>0 \qquad\mbox{for all } t\ge0.
\end{equation}
So the function $\bar S(\cdot)$ as defined in
(\ref{eq:fluid-accumulate}) has an inverse on the interval
$[0,\infty)$, which is denoted by $\bar T(\cdot)$. By the inverse
function theorem,
\[
\bar T'(v)=\bar Z(\bar T(v)) \qquad\mbox{for all }v\ge0.
\]
According to (\ref{eq:fluid-dym-eqn-S}) in Definition \ref{def:fluid},
we have
\[
\bar Z(t)=\mu\bigl(A_{\bar S(t)}\bigr)+\int_{0}^{t}[1-F(\bar S(s,t))]\,d\bar B(s).
\]
Perform the change of variables $u=\bar S(t)$ and $v=\bar S(s)$ to get
\begin{eqnarray*}
\bar Z(\bar T(u))&=&\mu(A_{u})+
\lambda\int_0^u[1-F(u-v)]\bar Z(\bar T(v))\,dv\\
&&{} -\int_{0}^{u}[1-F(u-v)]\,d\bar Q(\bar T(v)).
\end{eqnarray*}
Note that the function $\bar Q(\bar T(\cdot))$ has bounded variation
since it is the difference of two nondecreasing function $\bar B(\bar
T(\cdot))$ and $\lambda\bar T(\cdot)$. According to the integration
by parts formula provided by Lemma \ref{lem:LS-int-by-parts}, we
obtain
%
%
\begin{eqnarray}\label{eq:key-pre}\qquad\quad
\bar Z(\bar T(u))&=&\mu(A_{u})+
\lambda\beta\int_0^u\bar Z\bigl(\bar T(u-v)\bigr)\,dF_e(v)
-[1-F(0)]\bar Q(\bar T(u))\nonumber\\[-8pt]\\[-8pt]
&&{} +[1-F(u)]\bar Q(0) +\int_{0}^{u}\bar Q\bigl(\bar T(u-v)\bigr)\,dF(v),
\nonumber
\end{eqnarray}
where $F_e$ is the equilibrium distribution of $F$ which can be
written as $F_e(x)=\frac{1}{\beta}\int_0^x[1-F(y)]\,dy$. It follows
from (\ref{eq:critical-load-rho}) and (\ref{eq:cond-nu-11}) that
$\rho=1$ and $F(0)=0$, so we obtain the following key relationship:
%
%
\begin{eqnarray}\label{eq:key}\quad
\bar Q(\bar T(u))+\bar Z(\bar T(u))
&=&\xi(A_u)+\mu(A_u)+\int_0^{u}\bar Q\bigl(\bar T(u-v)\bigr)
\,dF(v)\nonumber\\[-8pt]\\[-8pt]
&&{} +\int_0^{u}\bar Z\bigl(\bar T(u-v)\bigr)\,dF_e(v)
\nonumber
\end{eqnarray}
for all $0\le u< \infty$. To simplify notation, denote
%
%
\begin{eqnarray}\label{eq:def-h}
h_{\xi,\mu}(u)&=&\xi(A_u)+\mu(A_u),\nonumber\\[-8pt]\\[-8pt]
x(u)&=&q(u)+z(u),\nonumber
\end{eqnarray}
where
%
%
\begin{eqnarray}
\label{eq:q-time-change}
q(u)&=&\bar Q(\bar T(u)),\\
\label{eq:z-time-change}
z(u)&=&\bar Z(\bar T(u)).
\end{eqnarray}
By (\ref{constr:policy-B}) and (\ref{constr:policy-S}), the above
equation can be written as
%
%
\begin{eqnarray}\label{eq:key-sim}
x(u)&=&h_{\xi,\mu}(u)+\int_0^{u}\bigl(x(u-v)-K\bigr)^+ \,dF(v)\nonumber\\[-8pt]\\[-8pt]
&&{}+\int_0^{u}\bigl(x(u-v)\wedge K\bigr)\,dF_e(v).\nonumber
\end{eqnarray}
By Lemma A.1 in \cite{ZDZ2009}, for any valid initial condition
$(\xi,\mu)$, the above integral equation has a unique solution which
is a c\`adl\`ag function. Our analysis on the limiting behavior of fluid
model solutions will be mainly based on (\ref{eq:key-sim}).

For a nonzero valid initial condition $(\xi,\mu)$ and an
$\varepsilon\in(0,1)$, define the $\varepsilon$-perturbation of it by
\[
(\xi_\varepsilon,\mu_\varepsilon)=\cases{
\biggl(\xi+\biggl(\varepsilon-\dfrac{K-\langle{1},{\mu}\rangle}{\langle{1},{\mu
}\rangle}\biggr)^+\mu,
\biggl(1+\dfrac{K-\langle{1},{\mu}\rangle}{\langle{1},{\mu}\rangle
}\wedge\varepsilon\biggr)\mu\biggr),\vspace*{2pt}\cr
\qquad\hspace*{84.7pt} \mbox{if $\langle{1},{\mu}\rangle<K$},\vspace*{2pt}\cr
(\xi+\varepsilon\mu,\mu),
\qquad\hspace*{26.74pt} \mbox{if $\langle{1},{\mu}\rangle=K,\xi=\mathbf{0}$},\vspace*{2pt}\cr
\bigl(\xi+\varepsilon(\xi+\mu),\mu\bigr),
\qquad \mbox{if $\xi\ne\mathbf{0}$},}
\]
and the $-\varepsilon$-perturbation of it by
\[
(\xi_{-\varepsilon},\mu_{-\varepsilon})=\cases{
\bigl(\xi,(1-\varepsilon)\mu\bigr),
&\quad if $\xi=\mathbf{0}$,\vspace*{2pt}\cr
\biggl(\biggl(1-\varepsilon\dfrac{\langle{\chi},{\xi+\mu}\rangle}{\langle
{\chi},{\xi}\rangle}\biggr)^+\xi,\cr
\qquad\biggl((1-\varepsilon)\dfrac{\langle{\chi},{\xi+\mu}\rangle}{\langle
{\chi},{\mu}\rangle}
1_{\{1<\varepsilon{\langle{\chi},{\xi+\mu}\rangle}/{\langle
{\chi},{\xi}\rangle}\}}\cr
\qquad\hspace*{78.7pt}+\,
1_{\{1\ge\varepsilon{\langle{\chi},{\xi+\mu}\rangle}/{\langle
{\chi},{\xi}\rangle}\}
}\biggr)\mu
\biggr), &\quad if $\xi\ne\mathbf{0}$.}
\]
The simple idea behind this complicated looking construction is that
(a) the perturbed state $(\xi_\varepsilon,\mu_\varepsilon)$ is still a valid
initial condition, (b) the workload of the perturbed state satisfies
%
%
\begin{equation}\label{eq:perturbe-workload}
\langle{\chi},{\xi_\varepsilon+\mu_\varepsilon}\rangle=(1+\varepsilon
)\langle{\chi},{\xi+\mu}\rangle\qquad
\mbox{for all }\varepsilon\in(-1,1),
\end{equation}
and (c) the function $h_{\xi_\varepsilon,\mu_\varepsilon}$ which is defined
based on $(\xi_\varepsilon,\mu_\varepsilon)$ in the same as (\ref{eq:def-h})
satisfies
%
%
\begin{equation}\label{eq:h-h_epsilon-bound}\quad
h_{\xi_{-\varepsilon},\mu_{-\varepsilon}}(u)
\le h_{\xi,\mu}(u)
\le h_{\xi_\varepsilon,\mu_\varepsilon}(u)\qquad
\mbox{for all }u\in\R_+,\varepsilon\in(0,1).
\end{equation}
The complication in the construction on the perturbation comes from
the requirement (\ref{eq:perturbe-workload}), which will provide
convenience in the proof of Lemmas \ref{lem:x-limit} and
\ref{lem:x-limit-unif}. Inequality (\ref{eq:h-h_epsilon-bound}) will
be used in the proof of the following lemma. Let
$x^{\varepsilon}(\cdot)$ denote the solution to (\ref{eq:key-sim}) with
$h_{\xi,\mu}$ replaced by $h_{\xi_\varepsilon,\mu_\varepsilon}$. We have
the following comparison.
\begin{lem}\label{lem:x-perturb}
Assume (\ref{eq:cond-nu-1}), (\ref{eq:cond-nu-11}) and
$(\xi,\mu)\ne(\mathbf{0},\mathbf{0})$. For all $\varepsilon\in(0,1)$,
\[
x^{-\varepsilon}(u)\le x(u)\le x^\varepsilon(u) \qquad\mbox{for all }u\ge0.
\]
\end{lem}
\begin{pf}
Let $u^*=\inf\{u\ge0\dvtx x(u)> x^\varepsilon(u)\}$. To prove $x(u)\le
x^\varepsilon(u)$, it is enough to show that $u^*=\infty$. Note that
$x(0)=h_{\xi,\mu}(0)<
h_{\xi_\varepsilon,\mu_\varepsilon}(0)=x^\varepsilon(0)$. By
right continuity of $x(\cdot)$ and $x^\varepsilon(\cdot)$, $u^*>0$.
Now, suppose\vadjust{\goodbreak} $u^*<\infty$. By (\ref{eq:key-sim}) and
(\ref{eq:h-h_epsilon-bound}) we have following bound estimation:
%
%
\begin{eqnarray}\label{eq:difference}
&& x^\varepsilon(u^*)-x(u^*)\nonumber\\
&&\qquad \ge\int_0^{u^*}\bigl[\bigl(x^\varepsilon(u^*-v)-K\bigr)^+-\bigl(x(u^*-v)-K\bigr)^+\bigr]\,d F(v)\\
&&\qquad\quad{} +\int_0^{u^*}\bigl[\bigl(x^\varepsilon(u^*-v)\wedge K\bigr)-\bigl(x(u^*-v)\wedge
K\bigr)\bigr]\,d F_e(v).\nonumber
\end{eqnarray}
Assumption (\ref{eq:cond-nu-11}) implies that $F(0)<1$. So there
exists $u'\in(0,u^*)$ such that
\[
\int_{u'-\delta}^{u'}dF_e(v)>0
\]
for all $0<\delta<u'$. By the definition of $u^*$, we have that
\[
\kappa=x^\varepsilon(u^*-u')-x(u^*-u')>0.
\]
By right continuity of $x(\cdot)$ and $x^\varepsilon(\cdot)$, we can
choose $\delta$ small enough such that
\[
x^\varepsilon(v)-x(v)\ge\frac{\kappa}{2} \qquad\mbox{for all }
u\in[u^*-u',u^*-u'+\delta].
\]
So by (\ref{eq:difference}), we have
\[
x^\varepsilon(u^*)-x(u^*)\ge\frac{\kappa}{2}
\int_{u'-\delta}^{u'}dF_e(v)>0.
\]
This contradicts the definition of $u^*$. So we must have that
$u^*=\infty$. The proof for the other inequality is completely
analogous.
\end{pf}

For the solution $x(\cdot)$ to (\ref{eq:key-sim}) with initial
condition $(\xi,\mu)$, define
%
%
\begin{equation}\label{eq:x-infinity}
x(\infty)=\frac{1}{\beta}(w-K\beta_e)^++\frac{1}{\beta_e}(w\wedge
K\beta_e),
\end{equation}
where $w=\langle{\chi},{\xi+\mu}\rangle$. The quantity $x(\infty
)$ can be
interpreted as the fluid system size corresponding to the equilibrium
state with workload $w$. We now use the above lemma and the key
renewal theorem to show the following convergence. To help with the
proof, we introduce the renewal function
\[
U(x)=\sum_{n=0}^\infty F^{*n}(x),
\]
associated with the distribution function $F$ (see Section V.2 in
\cite{Asmussen2003} for detailed discussion).
%
%
\begin{lem}\label{lem:x-limit}
Assume (\ref{eq:critical-load-rho}) and
(\ref{eq:cond-nu-1})--(\ref{eq:cond-nu-nonlattice}). The solution
$x(\cdot)$ to (\ref{eq:key-sim}) with initial condition
$(\xi,\mu)\in\init$ and $\langle{\chi},{\xi+\mu}\rangle<\infty
$ satisfies
\[
x(u)\to x(\infty) \qquad\mbox{as }u\to\infty.
\]
\end{lem}
\begin{pf}
We first study the case where $w=\langle{\chi},{\xi+\mu}\rangle
>K\beta_e$.
Convolve both sides of (\ref{eq:key-sim}) with, the renewal
function $U(\cdot)$ of $F(\cdot)$, to get
\[
x*U(u)=h_{\xi,\mu}*U(u)+(x-K)^+*F*U(u)+(x\wedge K)*F_e*U(u).
\]
Since $x=(x-K)^++(x\wedge K)$, by moving all terms containing
$(x-K)^+$ to the left and all terms containing $(x\wedge K)$ to the
right, we obtain
\[
\bigl(x(u)-K\bigr)^+*(1-F)*U(u)=h_{\xi,\mu}*U(u)-(x\wedge K)*(1-F_e)*U(u).
\]
This gives
%
%
\begin{eqnarray}\label{eq:tech-conv1}
\bigl(x(u)-K\bigr)^+&=&h_{\xi,\mu}*U(u)-K(1-F_e)*U(u)\nonumber\\[-8pt]\\[-8pt]
&&{} +[K-(x\wedge K)]*(1-F_e)*U(u).
\nonumber
\end{eqnarray}
Both $h_{\xi,\mu}(\cdot)$ and $1-F_e(\cdot)$ are directly Riemann
integrable since they are nonincreasing and integrable functions.
By the key renewal theorem, we have the convergence of the first two
terms on the right-hand side of (\ref{eq:tech-conv1}):
\begin{eqnarray*}
\lim_{u\to\infty}h_{\xi,\mu}*U(u)&=&\frac{w}{\beta},\\
\lim_{u\to\infty}K(1-F_e)*U(u)&=&\frac{K\beta_e}{\beta}.
\end{eqnarray*}
Note that $\frac{w}{\beta}-\frac{K\beta_e}{\beta}>0$ in this case,
and the last term in (\ref{eq:tech-conv1}) is always nonnegative.
So there exists $u_1>0$ such that
\[
\bigl(x(u)-K\bigr)^+>0 \qquad\mbox{for all }u\ge u_1.
\]
Equivalently, this
means that $K-(x(u)\wedge K)=0$ for all $u\ge u_1$. So the last
term in (\ref{eq:tech-conv1}) is nonnegative and can be bounded
above by
\[
\int_{u-u_1}^uK\,d[(1-F_e)*U(v)]=K[(1-F_e)*U(u)-(1-F_e)*U(u-u_1)],
\]
which converges to 0 by the key renewal theorem. So in this case we
have $\lim_{t\to\infty}x(u)=K+\frac{(w-K\beta_e)}{\beta}=x(\infty)$.

In the case where $w<K\beta_e$, we convolve both sides of
(\ref{eq:key-sim}) with $U_e(\cdot)$, the renewal function of
$F_e(\cdot)$ to get
\[
x*U_e(u)=h_{\xi,\mu}*U_e(u)+(x-K)^+*F*U_e(u)+(x\wedge K)*F_e*U_e(u).
\]
Again, since $x=(x-K)^++(x\wedge K)$, by moving all terms containing
$(x-K)^+$ to the right and all terms containing $(x\wedge K)$ to the
left, we obtain
%
%
\begin{equation}\label{eq:tech-conv2}
\bigl(x(u)\wedge K\bigr)=h_{\xi,\mu}*U_e(u)-(x-K)^+*(1-F)*U_e(u).
\end{equation}
By the key renewal theorem, the first term in the above converges:
\[
\lim_{u\to\infty}h_{\xi,\mu}*U_e(u)=\frac{w}{\beta_e}.
\]
Note that $\frac{w}{\beta_e}-K<0$ in this case, and the last term in
(\ref{eq:tech-conv2}) is always nonpositive. So there exists
$u_2>0$ such that
\[
\bigl(x(u)\wedge K\bigr)<K\qquad \mbox{for all } u\ge u_2.
\]
Equivalently, this means that $(x(u)-K)^+=0$ for all $u\ge u_2$.
The last term in (\ref{eq:tech-conv2}) is nonnegative and,
according to the upper bound (\ref{eq:upper-bound}), can be bounded
above by
\[
\int_{u-u_2}^u\frac{w}{\beta}\,d[(1-F)*U_e(v)]
=\frac{w}{\beta}[(1-F)*U_e(u)-(1-F)*U_e(u-u_2)],
\]
which converges to 0 by the key renewal theorem. So in this case, we
have $\lim_{t\to\infty}x(u)=\frac{w}{\beta_e}=x(\infty)$.

Now it only remains to study the case where $w=K\beta_e$. For any
$\varepsilon\in(0,1)$, let $(\xi_\varepsilon,\mu_\varepsilon)$ denote the
$\varepsilon$-perturbation of the initial condition $(\xi,\varepsilon)$,
introduced before Lemma \ref{lem:x-perturb}. It follows from
(\ref{eq:perturbe-workload}) that
\[
w_\varepsilon=\langle{\chi},{\xi_\varepsilon+\mu_\varepsilon}\rangle
=(1+\varepsilon
)\beta_eK.
\]
Following from the discussion of our first case:
\[
\lim_{u\to\infty}x^\varepsilon(u)=K+\frac{\varepsilon\beta_eK}{\beta}.
\]
By Lemma \ref{lem:x-perturb}, $x(u)\le x^\varepsilon(u)$ for all
$u\ge0$. So for all $\varepsilon>0$ there exists $u'_1$ such that when
$u\ge u'_1$
\[
x(u)\le K+\frac{\varepsilon\beta_eK}{\beta}+\varepsilon
= K+\biggl(\frac{\beta_eK}{\beta}+1\biggr)\varepsilon.
\]
Similarly, we have the $-\varepsilon$-perturbation
$(\xi_{-\varepsilon},\mu_{-\varepsilon})$ for $-\varepsilon\in(-1,0)$. It
follows from (\ref{eq:perturbe-workload}) that
\[
w_{-\varepsilon}=\langle{\chi},{\xi_{-\varepsilon}+\mu_{-\varepsilon
}}\rangle
=(1-\varepsilon)\beta_eK.
\]
Following from the discussion of our first case,
\[
\lim_{u\to\infty}x^{-\varepsilon}(u)=K-\varepsilon K.
\]
By Lemma \ref{lem:x-perturb}, $x(u)\ge x^{-\varepsilon}(u)$ for all
$u\ge0$. So for all $\varepsilon>0$ there exists $u'_2$ such that when
$u\ge u'_2$
\[
x(u)\ge K -\varepsilon K-\varepsilon=K-(K+1)\varepsilon.
\]
Summarizing this case, we have $\lim_{t\to\infty}x(u)=K=x(\infty)$.
\end{pf}
%
%
\begin{lem}\label{lem:Ax-converges}
Assume (\ref{eq:critical-load-rho}) and
(\ref{eq:cond-nu-1})--(\ref{eq:cond-nu-nonlattice}). Let
$(\bar\buf(\cdot),\bar\ser(\cdot))$ be the solution to the fluid
model $(K,\lambda,\nu)$ with initial condition $(\xi,\mu)\in\init$
and $\langle{\chi},{\xi+\mu}\rangle<\infty$. Let $w=\langle{\chi
},{\mu}\rangle$. We have
as $t\to\infty$,
%
%
\begin{eqnarray}
\label{eq:pov-B-conv}
\pov\biggl[\bar\buf(t),\frac{(w-K\beta_e)^+}{\beta}\nu\biggr]&\to&0,\\
\label{eq:Ax-S-conv}
\sup_{y\in(0,\infty)}\biggl|\bar\ser(t)(A_y) - \frac{w\wedge
K\beta_e}{\beta_e}\nu_e(A_y)\biggr|&\to&0.
\end{eqnarray}
\end{lem}
\begin{pf}
If $w=0$, the result holds trivially. Now assume that $w\ne0$. Let
\begin{eqnarray*}
q(\infty)&=&\bigl(x(\infty)-K\bigr)^+=\frac{(w-K\beta_e)^+}{\beta},\\
z(\infty)&=&x(\infty)\wedge K=\frac{w\wedge K\beta_e}{\beta_e},
\end{eqnarray*}
where $x(\infty)$ is defined in (\ref{eq:x-infinity}). Based on the
fluid dynamic equation (\ref{eq:fluid-dym-eqn-B}) and the fact that
$\nu$ is a probability measure, we have
\[
|\bar\buf(t)(A)-q(\infty)\nu(A)|
=|[\bar Q(t)-q(\infty)]\nu(A)|\le|\bar Q(t)-q(\infty)|,
\]
for all Borel set $A$. This implies that
\[
\bar\buf(t)(A)\le q(\infty)\nu(A)+ |\bar Q(t)-q(\infty)|
\le q(\infty)\nu\bigl(A^{|\bar Q(t)-q(\infty)|}\bigr)
+|\bar Q(t)-q(\infty)|,
\]
where $A^\kappa$ denote the $\kappa$-enlargement as introduced
in Section \ref{subsec:notatioin}. Similarly, we have
\[
q(\infty)\nu(A)\le\bar\buf(t)\bigl(A^{|\bar Q(t)-q(\infty)|}\bigr)
+|\bar Q(t)-q(\infty)|.
\]
By the change of variable $u=\bar S(t)$ [$t=\bar T(u)$] and the
definition of the Prohorov metric,
\[
\pov[\bar\buf(t),q(\infty)\nu]\le|q(u)-q(\infty)|.
\]
By Lemma \ref{lem:x-limit}, there exists $u_1>0$ such that when
$u>u_1$ we have $|q(u)-q(\infty)|\le\varepsilon$. So for all
$\varepsilon>0$, there exists $t_1=Ku_1\ge\bar T(u_1)$ such that
%
%
\begin{equation}\label{eq:ax-buffer-conv}
\pov[\bar\buf(t),q(\infty)\nu]
\le\varepsilon\qquad\mbox{for all }t>t_1.
\end{equation}

It remains to study the limiting behavior of $\bar\ser(\cdot)$.
Perform the change of variable $u=\bar S(t)$ ($t=\bar T(u)$) to the
fluid dynamic equation (\ref{eq:fluid-dym-eqn-S}), we get
\[
\bar\ser(\bar T(u))(A_y)=\mu(A_y+u)
+\int_{0}^{u}\nu(A_y+u-v) \,d[\lambda\bar T(v)-q(v)].
\]
Due to the fact that $\rho=1$, we have
$\frac{1}{\beta}=\lambda$. Thus,
$z(\infty)\nu_e(A_y)=\lambda\int_{0}^{u}\nu(A_y+u-v)
z(\infty)\,dv+z(\infty)\nu_e(A_y+u)$. Since $d\bar T(v)=z(v)\,dv$, we
then have the following bound estimation:
%
%
\begin{eqnarray}\label{eq:tech-3-2}
&& |\bar\ser(\bar T(u))(A_y)-z(\infty)\nu_e(A_y)
|\nonumber\\
&&\qquad\le\mu(A_y+u)+\biggl|\int_{0}^{u}\nu(A_y+u-v) \,d
q(v)\biggr|\nonumber\\[-8pt]\\[-8pt]
&&\qquad\quad{} +\biggl|\lambda\int_{0}^{u}\nu(A_y+u-v) [z(\infty
)-z(v)]\,dv\biggr|\nonumber\\
&&\qquad\quad{}+z(\infty)\nu_e(A_y+u).\nonumber
\end{eqnarray}
%
It is clear that the first and the last terms on the right-hand side of
(\ref{eq:tech-3-2}) vanishes as $u\to\infty$.
Recall that $F$ is the distribution function corresponding to the
measure~$\nu$. By integration by parts (see
Lemma \ref{lem:LS-int-by-parts}), the second term on the right-hand
side of (\ref{eq:tech-3-2}) can be written as
\[
\biggl|[1-F(y)]q(u)-[1-F(y+u)]q(0)+\int_{0}^{u}q(v)\,dF(y+u-v)\biggr|,
\]
which is less than or equal to
\begin{eqnarray*}
&& |[1-F(y)]q(u)-[F(y+u)-F(y)]q(\infty)|
+|[1-F(y+u)]q(0)|\\
&&\qquad{} +\int_{0}^{u}|q(u-v)-q(\infty)|\,dF(y+v).
\end{eqnarray*}
By convergence of $x(\cdot)$, for all $\varepsilon>0$ there exists a
$u_1$ such that $|x(v)-x(\infty)|<\varepsilon$ if $v\ge u_1$. For all
$\varepsilon>0$, we can choose $u_2>0$ such that $1-F(u_2)<\varepsilon$.
When $u>u_1+u_2$, the above inequality can be further bounded by
\begin{eqnarray*}
&& [1+q(0)+q(\infty)]\varepsilon
+\int_{0}^{u-u_1}|q(u-v)-q(\infty)|\,dF(y+v)\\
&&\quad{} +\int_{u-u_1}^{u}|q(u-v)-q(\infty)|\,dF(y+v)\\
&&\qquad\le[1+q(0)+q(\infty)]\varepsilon+[F(u-u_1+y)-F(y)]\varepsilon\\
&&\qquad{} +{\sup_{v\ge0}}|q(v)-q(\infty)|[F(u)-F(u-u_1)]\\
&&\qquad\le\biggl(2+\frac{2w}{\beta}\biggr)\varepsilon,
\end{eqnarray*}
where the last inequality is due to (\ref{eq:upper-bound}).
When $u>u_1+u_2$, the third term in (\ref{eq:tech-3-2}) can be
written as
\[
\biggl|\lambda\int_{0}^{u_1}\nu(A_y+u-v) [z(\infty)-z(v)]\,dv
+\lambda\int_{u_1}^{u}\nu(A_y+u-v) [z(\infty)-z(v)]\,dv\biggr|,
\]
which is bounded above by
\begin{eqnarray*}
&&\lambda\sup_{0\le v\le u_1}|z(v)-z(\infty)|[1-F(u-u_1)]
+\varepsilon\lambda\int_{u_1}^u [1-F(v)]\,dv\\
&&\qquad\le\lambda K\varepsilon+\lambda\beta\varepsilon,
\end{eqnarray*}
where the last inequality is due to the bound $z(u)\le K$ for all $u\ge
0$. So for all $\varepsilon>0$ there exists a $t_2=K(u_1+u_2)\ge\bar
T(u_1+u_2)$ such that
%
%
\begin{equation}\label{eq:ax-server-conv}
{\sup_{y\in(0,\infty)}}|\bar\ser(t)(A_y)-z(\infty)\nu
_e(A_y)|<\varepsilon\qquad
\mbox{for all }t\ge t_2.
\end{equation}
\upqed\end{pf}
\begin{pf*}{Proof of Theorem \ref{thm:convergence-to-invariant-mani},
part I}
Note that $\langle{\chi},{\frac{w\wedge K\beta_e}{\beta_e}\nu
_e}\rangle\le
K\beta_e<\infty$. By the workload conversing property,
$\langle{\chi},{\bar Z(\cdot)}\rangle\le w<\infty$. According to
Lemma \ref{lem:conv-proh}, (\ref{eq:Ax-S-conv}) in
Lemma \ref{lem:Ax-converges} implies that
\[
\pov\biggl[\bar\ser(t),\frac{w\wedge K\beta_e}{\beta_e}\nu_e\biggr]\to0
\qquad\mbox{as }t\to\infty.
\]
This and (\ref{eq:pov-B-conv}) implies the convergence result in
Theorem \ref{thm:invariant-manifold}.
\end{pf*}

\subsection{Uniform convergence to equilibrium states}
\label{subsec:unif-conv-invar-manif}
The convergence in the previous subsection depends on the initial
condition $(\xi,\mu)$. We now show that the convergence is uniform
for all initial conditions in the set $\init_M^p$ defined in
Theorem~\ref{thm:convergence-to-invariant-mani}.

To emphasize the dependency on the initial condition, we use
$\solx(\init_M^p)$ to denote the set of solutions to equation
(\ref{eq:key-sim}) with input function $h_{\xi,\mu}$ induced by
initial condition $(\xi,\mu)\in\init_M^p$, and $\sol(\init_M^p)$ to
denote the set of solutions to the fluid model $(K,\lambda,\nu)$ with
initial condition $(\xi,\mu)\in\init_M^p$.
%
%
\begin{lem}\label{lem:x-limit-unif}
Assume (\ref{eq:critical-load-rho}) and
(\ref{eq:cond-nu-1})--(\ref{eq:cond-nu-nonlattice}). For each
$\varepsilon>0$ there exists an $l^*>0$ such that when $u\ge l^*$,
\[
{\sup_{x(\cdot)\in\solx(\init_M^p)}}|x(u)-x(\infty)|<\varepsilon.
\]
\end{lem}
\begin{pf}
To prove this lemma, we need to adjust the proof of
Lemma \ref{lem:x-limit} with the assistance of
Lemma \ref{lem:key-renewal-thm-unif}.

Let $\mathscr H_M=\{h_{\xi,\mu}\dvtx(\xi,\mu)\in\init_M^p\}$. By the
definition of the set $\init_M^p$ in
Theorem~\ref{thm:convergence-to-invariant-mani}, $\mathscr H_M$ is
the set of nonincreasing functions satisfying condition
(\ref{eq:bound-h-0}) and (\ref{eq:unif-int}) in
Lemma \ref{lem:key-renewal-thm-unif}. For any $\varepsilon>0$, divide
the set $\init_M^p$ into three parts,
\[
\init_M^p=\init^+_\varepsilon\cup\init^0_\varepsilon\cup\init
^-_\varepsilon,
\]
where
\begin{eqnarray*}
\init^+_\varepsilon&=&\{(\xi,\mu)\in\init_M^p\dvtx\langle{\chi},{\xi
+\mu}\rangle
\ge K\beta_e(1+\varepsilon)\},\\
\init^-_\varepsilon&=&\{(\xi,\mu)\in\init_M^p\dvtx\langle{\chi},{\xi
+\mu}\rangle
\le K\beta_e(1-\varepsilon)\},
\end{eqnarray*}
and
$\init^0_\varepsilon=\init_M^p\setminus(\init^+_\varepsilon\cup\init
^-_\varepsilon)$.

We first focus on the set $\init^+_\varepsilon$. By doing the same
algebra as in the proof of Lemma \ref{lem:x-limit}, we see that
(\ref{eq:tech-conv1}) holds for any $(\xi,\mu)\in\init^+_\varepsilon$.
By Lemma \ref{lem:key-renewal-thm-unif} and the key renewal theorem,
there exists a $u^*_1$ such that
\begin{eqnarray*}
\sup_{(\xi,\mu)\in\init^+_\varepsilon}
\biggl|h_{\xi,\mu}*U(u)-\frac{\langle{\chi},{\xi+\mu}\rangle}{\beta}\biggr|
&<&\frac{K\beta_e}{4\beta}\varepsilon,\\
\biggl|K(1-F_e)*U(u)-K\frac{\beta_e}{\beta}\biggr|
&<&\frac{K\beta_e}{4\beta}\varepsilon,
\end{eqnarray*}
for all $u\ge u^*_1$. So for the first two terms on the right-hand
side of (\ref{eq:tech-conv1}), we have
\[
h_{\xi,\mu}*U(u)-K(1-F_e)*U(u)\ge
\frac{K\beta_e(1+\varepsilon)}{\beta}
-K\frac{\beta_e}{\beta}-\frac{K\beta_e}{2\beta}\varepsilon>0,
\]
for all $(\xi,\mu)\in\init^+_\varepsilon$ and $u>u^*_1$. Note that the
last term in (\ref{eq:tech-conv1}) is always nonnegative. So when
$u\ge u^*_1$ we have $(x(u)-K)^+>0$ [or equivalently $K-(x(u)\wedge
K)=0$], for all $(\xi,\mu)\in\init^+_\varepsilon$. So the last term on
the right-hand side of (\ref{eq:tech-conv1}) is nonnegative and can
be bounded above by
\[
\int_{u-u^*_1}^uK\,d[(1-F_e)*U(v)]=K[(1-F_e)*U(u)-(1-F_e)*U(u-u^*_1)],
\]
which converges to 0 as $u\to\infty$ by the key renewal theorem. So
there exists a \mbox{$u'_1>0$} such that when $u>u'_1$, the absolute value
of third term on the right-hand side of (\ref{eq:tech-conv1}) is
bounded by $\frac{K\beta_e}{2\beta}\varepsilon$. Let
$l^*_1=\max(u'_1,u^*_1)$. By (\ref{eq:tech-conv1}) and summarizing
the above, we obtain
\[
{\sup_{(\xi,\mu)\in\init^+_\varepsilon}}|x(u)-x(\infty)|
<\frac{K\beta_e}{\beta}\varepsilon\qquad\mbox{for all }u>l^*_1.
\]

Next, we consider the set $\init^-_\varepsilon$. By doing the same
algebra as in the proof of Lemma \ref{lem:x-limit}, we see that
(\ref{eq:tech-conv2}) holds for any $(\xi,\mu)\in\init^-_\varepsilon$.
By Lemma \ref{lem:key-renewal-thm-unif}, there exists a $u^*_2$ such
that
\[
\sup_{(\xi,\mu)\in\init^-_\varepsilon}\biggl|h_{\xi,\mu}*U_e(u)
-\frac{\langle{\chi},{\xi+\mu}\rangle}{\beta_e}\biggr|<\frac
{K}{2}\varepsilon
\]
for all $u>u^*_2$. So we have
\[
h_{\xi,\mu}*U_e(u)\le\frac{K\beta_e(1-\varepsilon)}{\beta_e}
+\frac{K}{2}\varepsilon<K,
\]
for all $(\xi,\mu)\in\init^-_\varepsilon$ and $u>u^*_2$. Note that the
last term in (\ref{eq:tech-conv2}) is always nonpositive. So when
$u\ge u^*_2$ we have $x(u)< K$ [or equivalently, $(x(u)-K)^+=0$] for
all $(\xi,\mu)\in\init^-_\varepsilon$. So by (\ref{eq:upper-bound}),
the absolute value of the last term on the right-hand side of
(\ref{eq:tech-conv2}) can be bounded by
\[
\int_{u-u^*_2}^u\frac{w}{\beta}\,d[(1-F)*U_e(v)]
=\frac{w}{\beta}[(1-F)*U_e(u)-(1-F)*U_e(u-u^*_2)],
\]
which converges to 0 as $u\to\infty$ by the key renewal theorem. So
there exists a $u'_2>0$ such that when $u>u'_2$, the last term on
the right-hand side of (\ref{eq:tech-conv2}) is bounded by
$\frac{K}{2}\varepsilon$. Let $l^*_2=\max(u'_2,u^*_2)$. By
(\ref{eq:tech-conv2}) and summarizing the above,
\[
{\sup_{(\xi,\mu)\in\init^-_\varepsilon}}|x(u)-x(\infty)|<K\varepsilon
\qquad\mbox{for all }u>l^*_2.
\]

It only remains to deal with the set $\init^0_\varepsilon$. We can
restrict $\varepsilon<1/3$, since we are only interested in small ones.
According to (\ref{eq:perturbe-workload}), for any
$(\xi,\mu)\in\init^0_\varepsilon$, we have
$h_{\xi_{3\varepsilon},\mu_{3\varepsilon}}\in\init^+_\varepsilon$ and
$h_{\xi_{-3\varepsilon},\mu_{-3\varepsilon}}\in\init^-_\varepsilon$. Denote
$x^+(\cdot)$ and $x^-(\cdot)$ the solutions to (\ref{eq:key-sim})
corresponding to $h_{\xi_{3\varepsilon},\mu_{3\varepsilon}}$ and
$h_{\xi_{-3\varepsilon},\mu_{-3\varepsilon}}$, respectively. By
Lemma \ref{lem:x-perturb},
\[
x^-(u)<x(u)<x^+(u) \qquad\mbox{for all }u\ge0.
\]
Note that in
this case, the workload $w\le(1+\varepsilon)\beta_eK\le2\beta_eK$. By
(\ref{eq:x-infinity}), we have that
\begin{eqnarray*}
x^+(\infty)
&\le& x(\infty)+\max\biggl(\frac{1}{\beta},\frac{1}{\beta_e}\biggr)3\varepsilon w
\le x(\infty)+\max\biggl(\frac{1}{\beta},\frac{1}{\beta_e}\biggr)6\beta
_eK\varepsilon
,\\
x^-(\infty)
&\ge& x(\infty)-\max\biggl(\frac{1}{\beta},\frac{1}{\beta_e}\biggr)3\varepsilon w
\ge x(\infty)-\max\biggl(\frac{1}{\beta},\frac{1}{\beta_e}\biggr)6\beta
_eK\varepsilon.
\end{eqnarray*}
According to the above two cases, when $u>l^*=\max(l^*_1,l^*_2)$,
\begin{eqnarray*}
x(u)&\le& x^+(\infty)+\frac{K\beta_e}{\beta}\varepsilon
\le x(\infty)
+\max\biggl(\frac{1}{\beta},\frac{1}{\beta_e}\biggr)6\beta_eK\varepsilon
+\frac{K\beta_e}{\beta}\varepsilon,\\
x(u)&\ge& x^-(\infty)-K\varepsilon
\ge x(\infty)
-\max\biggl(\frac{1}{\beta},\frac{1}{\beta_e}\biggr)6\beta_eK\varepsilon
-K\varepsilon,
\end{eqnarray*}
for all $(\xi,\mu)\in\init^0_\varepsilon$. This means that
\[
{\sup_{(\xi,\mu)\in\init^0_\varepsilon}}|x(u)-x(\infty)|<C\varepsilon
\qquad\mbox{for all }u>l^*,
\]
where $C=\max(\frac{1}{\beta},\frac{1}{\beta_e})6\beta_eK
+\frac{K\beta_e}{\beta}+K$.
\end{pf}
%
%
\begin{lem}\label{lem:Ax-converges-unif}
Assume (\ref{eq:critical-load-rho}) and
(\ref{eq:cond-nu-1})--(\ref{eq:cond-nu-nonlattice}). For all
$\varepsilon>0$ there exists an $L^*>0$ such that when $t\ge L^*$,
%
%
\begin{eqnarray}\label{eq:pov-B-conv-unif}
\sup_{(\bar\buf(\cdot),\bar\ser(\cdot))\in\sol(\init_M^p)}
\pov\biggl[\bar\buf(t),\frac{(w-K\beta_e)^+}{\beta}\nu\biggr]&<&\varepsilon,\\
\label{eq:Ax-S-conv-unif}
\sup_{(\bar\buf(\cdot),\bar\ser(\cdot))\in\sol(\init_M^p)}
\sup_{y\in(0,\infty)}\biggl|\bar\ser(t)(A_y)
-\frac{w\wedge K\beta_e}{\beta_e}\nu_e(A_y)\biggr|&<&\varepsilon.
\end{eqnarray}
\end{lem}
\begin{pf}
The proof of this corollary is almost the same as the proof of
Lem\-ma~\ref{lem:Ax-converges}.
Just note that by Lemma \ref{lem:x-limit-unif}, the $t_1$ in
(\ref{eq:ax-buffer-conv}) and the $t_2$ in (\ref{eq:ax-server-conv})
are good for all $(\xi,\mu)\in\init_M^p$. With $L^*=\max(t_1,t_2)$,
the result of this lemma immediately follows.
\end{pf}
\begin{pf*}{Proof of Theorem \ref{thm:convergence-to-invariant-mani}, part
II}
Now we use Lemma \ref{lem:Ax-converges-unif} to show the uniform
convergence result. Note that $\langle{\chi},{\frac{w\wedge K\beta
_e}{\beta_e}\nu_e}\rangle\le K\beta_e<\infty$. By the definition
of $\init^p_M$, for any
$(\bar\buf(\cdot),\bar\ser(\cdot))\in\sol(\init_M^p)$,
$\langle{\chi},{\bar Z(\cdot)}\rangle<M<\infty$. According to
Lem\-ma~\ref{lem:conv-proh}, (\ref{eq:Ax-S-conv-unif}) in
Lemma \ref{lem:Ax-converges-unif} implies that for all $\varepsilon>0$
there exists an $L^*_1$ such that when $t\ge L^*_1$,
\[
\sup_{(\bar\buf(\cdot),\bar\ser(\cdot))\in\sol(\init_M^p)}
\pov\biggl[\bar\ser(t),\frac{w\wedge K\beta_e}{\beta_e}\nu_e\biggr]<\varepsilon.
\]
The uniform convergence follows from the above and
(\ref{eq:pov-B-conv-unif}).
\end{pf*}

\section{Shifted fluid scaling and precompactness}
\label{sec:tightness}

The objective of this section is to show the precompactness property,
Theorem \ref{thm:precompactness} at the end of this section, for the
sequence of shifted fluid scaled processes,
which is defined in the following
section.

\subsection{Shifted fluid scaling}\label{sec:shift-fluid-scal}

Much of our understanding of the diffusion scaled process will be
derived from results about the \textit{shifted fluid scaled process},
which is defined by
%
%
\begin{equation}
\frm\buf(t)=\frac{1}{r}\buf^r(rm+rt),\qquad
\frm\ser(t)=\frac{1}{r}\ser^r(rm+rt),
\end{equation}
for all $m\in\N$ and $t\ge0$. To see the relationship between these
two scalings, consider the diffusion scaled process on the interval
$[0,T]$, which corresponds to the interval $[0,r^2T]$ for the unscaled
process. Fix a constant $L>1$, the interval will be covered by the
$\lfloor{rT}\rfloor+1$ overlapping intervals
\[
[rm,rm+rL],\qquad m=0,1,\ldots,\lfloor{rT}\rfloor.
\]
For each
$t\in[0,T]$, there exists an $m\in\{0,\ldots,\lfloor{rT}\rfloor\}$ and
$s\in[0,L]$ (which may not be unique) such that $r^2t=rm+rs$. Thus,
%
%
\begin{equation}\label{eq:relation-frm-dr}
\dr{\buf}(t)=\frm{\buf}(s),\qquad \dr{\ser}(t)=\frm{\ser}(s).
\end{equation}
This will serve as a key relationship between fluid and diffusion
scaled processes.

We are also interested in shifted fluid scaled versions of other
processes, like the workload and system size processes. Note that
$Q^r(\cdot)$, $Z^r(\cdot)$, $X^r(\cdot)$, $W^r(\cdot)$ and
$S^r(\cdot,\cdot)$ are actually functions of
$(\buf^r(\cdot),\ser^r(\cdot))$, so the scaling for these quantities
is defined as the functions of the corresponding scaling for
$(\buf^r(\cdot),\ser^r(\cdot))$, that is,
%
%
\begin{eqnarray}
\frm{Q}(t)&=&\langle{1},{\frm\buf(t)}\rangle=\frac
{1}{r}Q^r(rm+rt),\\
\frm{Z}(t)&=&\langle{1},{\frm\ser(t)}\rangle=\frac
{1}{r}Z^r(rm+rt),\\
\frm{X}(t)&=&\langle{1},{\frm\ser(t)+\frm\ser(t)}\rangle=\frac
{1}{r}X^r(rm+rt),\\
\frm{W}(t)&=&\langle{\chi},{\frm\buf(t)+\frm\ser(t)}\rangle
=\frac{1}{r}W^r(rm+rt),
\end{eqnarray}
for all $0\le s\le t$. We define the shifted fluid scaling for the
arrival process as
%
%
\begin{equation}
\frm E(t)=\frac{1}{r}E^r(rm+rt),
\end{equation}
for all $t\ge0$. By (\ref{eq:B(t)}), the shifted fluid scaling for
$B^r(\cdot)$ is
%
%
\begin{equation}
\frm B(t)=\frm E(t)-\frm Q(t),
\end{equation}
for all $t\ge0$. To shorten the notation, for all $0\le s\le t$,
denote
%
%
\begin{equation}\label{eq:shifted-B}\qquad
\frm E(s,t)=\frm E(t)-\frm E(s),\qquad
\frm B(s,t)=\frm B(t)-\frm B(s).
\end{equation}

A shifted fluid scaled version of the stochastic dynamic equations
(\ref{eq:stoc-dym-eqn-B}) and (\ref{eq:stoc-dym-eqn-S}) can be written
as, for $0\le s\le t$,
%
%
\begin{eqnarray}
\label{eq:mechanism-B}
\frm\buf(t)(A)&=&\frm\buf(s)(A)
+\frac{1}{r}\sum_{i=r\frm E(s)+1}^{r\frm
E(t)}\delta_{v^r_i}(A)\nonumber\\[-8pt]\\[-8pt]
&&{} -\frac{1}{r}\sum_{i=r\frm B(s)+1}^{r\frm B(t)}
\delta_{v^r_i}(A),\qquad
A\subseteq[0,\infty),\nonumber\\
\label{eq:mechanism-S}
\frm\ser(t)(A)&=&\frm\ser(s)\bigl(A+S^r(rm+rs,rm+rt)\bigr)\nonumber\\
&&{} +\frac{1}{r}\sum_{i=r\frm B(s)+1}^{r\frm B(t)}
\delta_{v^r_i}\bigl(A+S^r(\tau^r_i,rm+rt)\bigr),\\
&&\eqntext{A\subseteq(0,\infty).}
\end{eqnarray}
Please note that $\frm\ser(t)(\{0\})=0$ for all $t\ge0$ according to
our definition. The dynamics of the system is determined by the above
equations. Equation (\ref{eq:mechanism-B}) says that the status of
the buffer at time $t$ equals the status at time $s$ plus what has
arrived to the buffer and minus what has left from the buffer during
time interval $(s,t]$. Those jobs who left buffer enter service; the
service process has been taken care of by shifting the set $A$ by the
cumulative service amount $S^r(\tau_i,rm+rt)$ that the $i$th job
receives. This corresponds to the second term on the right-hand side
of (\ref{eq:mechanism-S}). This plus the status at time $s$ shifted
by accumulative service amount $S^r(rm+rs,rm+rt)$ is equal to the
status of the server at time $t$, as indicated
in~(\ref{eq:mechanism-S}).

\subsection{Preliminary estimates}\label{subsec:some-bound-estimates}

We first establish some bounds which will be useful for later
discussion. The following lemma gives some bound on the arrival
processes.
%
%
\begin{lem}\label{lem:arrival-process}
Assume (\ref{eq:cond-A}) and (\ref{eq:cond-lambda}). Fix $T>0$ and
$L>1$. For all $\varepsilon$, \mbox{$\varepsilon'>0$}, there exists an $r_0$ such
that whenever $r\ge r_0$,
%
%
\begin{equation}\label{ineq:arrival-bound}
\mathbb{P}^r\Bigl({{\max_{m\le\lfloor{rT}\rfloor}\sup_{s,t\in
[0,L]}}|E^{r,m}(s,t)-\lambda(t-s)|>\varepsilon' }\Bigr)<\varepsilon.
\end{equation}
\end{lem}
\begin{pf}
Let $t'=\frac{m+t}{r}$ and $s'=\frac{m+s}{r}$. Note that
$\max_{m\le\lfloor{rT}\rfloor}\sup_{t\in[0,L]}\frac{m+t}{r}<T+1$
for all large
$r$, and $0\le s,t\le L$ is the same as $|t'-s'|\le L/r$. For any
$\delta>0$, there exists an $r'_0$ such that $L/r<\delta$ for all
$r\ge r_0$, so the left-hand side of (\ref{ineq:arrival-bound}) can be
bounded above by
%
%
\begin{equation}\label{eq:tech-arrival-bound}
\mathbb{P}^r\biggl({\sup_{s',t'\in[0,T+1],|s'-t'|<\delta} \biggl|\frac
{1}{r} E^r(r^2t')-\lambda rt' -\biggl(\frac{1}{r}E^r(r^2s')-\lambda
rs'\biggr)\biggr|>\varepsilon' }\biggr)\hspace*{-35pt}
\end{equation}
for all $r\ge r'_0$. By the assumptions (\ref{eq:cond-A}) and
(\ref{eq:cond-lambda}) on the arrival process,
$\{\frac{1}{r}E(r^2\cdot)-\lambda r\cdot\}$ converges in
distribution to the Brownian motion $E^*(\cdot)$. Since a Brownian
motion is almost surely continuous, we conclude that
(\ref{eq:tech-arrival-bound}) converges to zero as $\delta\to0$.
Then the inequality (\ref{ineq:arrival-bound}) follows immediately.
\end{pf}

Here is a remark that will facilitate some arguments later on. The
$\varepsilon'$ and $\varepsilon$ in (\ref{ineq:arrival-bound}) can be
replaced by $\varepsilon_E(r)$, which is a function of $r$ that vanishes
at infinity. Here is the proof. For each index $r$ let
\[
H_r=\{\delta>0\dvtx\mbox{(\ref{ineq:arrival-bound})  is true for }
\varepsilon'=\varepsilon=\delta\}.
\]
Clearly, $H_r$ is not empty since $2\in H_r$. Let $\varepsilon_E(r)=\inf
H_r$ for each $r\ge0$. Assume that $\varepsilon_E(r)$ does not vanish
at infinity. There exists a $\delta>0$ and a subsequence
$\{r_n\}_{n=1}^\infty$ which increases to infinity such that
%
%
\begin{equation}\label{eq:epsilon-r}
\varepsilon_E(r_n)>\delta\qquad\mbox{for all }n\ge0.
\end{equation}
However, by Lemma \ref{lem:arrival-process}, for
$\varepsilon'=\varepsilon=\delta/2$ there exists an $r_\delta$ such that
when $r_n\ge r_\delta$, (\ref{ineq:arrival-bound}) must hold. This
contradicts (\ref{eq:epsilon-r}). Based on this, we construct
%
%
\begin{equation}\label{eq:omegaE-def}
\Omega^r_E=\Bigl\{
{\max_{m\le\lfloor{rT}\rfloor}\sup_{s,t\in[0,L]}}
|E^{r,m}(s,t)-\lambda(t-s)|\le\varepsilon_E(r)
\Bigr\}.
\end{equation}
According to Lemma \ref{lem:arrival-process}, we have that
%
%
\begin{equation}\label{eq:omegaE}
\lim_{r\to\infty}\mathbb{P}^r({\Omega^r_E})=1.
\end{equation}

Recall the Glivenko--Cantelli estimate in
Lemma \ref{lem:glivenko-cantelli}. By the same argument as in the
above, for fixed constant $M_1,L_1$, there exists a function
$\varepsilon_{\mathrm{GC}}(\cdot)$, which vanishes at infinity, such that
the probability inequality in Lemma \ref{lem:glivenko-cantelli} holds
with $\varepsilon$ and $\varepsilon'$ replaced by this function. In other
words, if we denote
%
%
\begin{eqnarray}\label{eq:omega-GC-def}\qquad
\Omega^r_{\mathrm{GC}}(M_1,L_1)&=&\Bigl\{
\max_{-rM_1< n< r^2M_1}\sup_{l\in[0,L_1]}\sup_{f\in\testfn}
|\langle{f},{\fr\eta(n,l)}\rangle-l\langle{f},{\nu^r}\rangle|\nonumber\\[-8pt]\\[-8pt]
&&\hspace*{182pt}\le
\varepsilon_{\mathrm{GC}}(r)
\Bigr\},\nonumber
\end{eqnarray}
where $\fr\eta(n,l)$ is defined in (\ref{eq:def-eta}) and $\testfn$ is
defined in (\ref{eq:set-of-test-fns}), then for any fixed constant
$M_1,L_1$,
%
%
\begin{equation}\label{eq:omega-GC}
\lim_{r\to\infty}\mathbb{P}^r({\Omega^r_{\mathrm
{GC}}(M_1,L_1)})=1.
\end{equation}

Now, we use the above result and Proposition \ref{prop:workload} to
obtain a bound on the queue length processes.
%
%
\begin{lem}\label{lem:Q-bound}
Assume (\ref{eq:cond-A})--(\ref{eq:cond-S-u}),
(\ref{eq:cond-initial}), (\ref{eq:cond-initial-u}) and
(\ref{eq:cond-HT}). Fix $T>0$ and $L>1$. For all $\eta>0$, there
exists a constant $M>0$ such that
\[
\liminf_{r\to\infty}\mathbb{P}^r\Bigl({\max_{m\le\lfloor
{rT}\rfloor}\sup
_{t\in[0,L]}\frm Q(t)\le M }\Bigr)\ge1-\eta.
\]
\end{lem}
\begin{pf}
Since $\frac{\lfloor{rT}\rfloor+L}{r}<T+1$ for all large enough $r$,
it is
enough to prove the following inequality:
\[
\liminf_{r\to\infty}\mathbb{P}^r\Bigl({\sup_{t\in[0,T+1]}\dr
Q(t)\le M}\Bigr)\ge
1-\eta.
\]
Suppose this is not true, then there exists an $\eta>0$ such that
for any $M$,
\[
\liminf_{r\to\infty}\mathbb{P}^r\Bigl({\sup_{t\in[0,T+1]}\dr
Q(t)>M}\Bigr)>\eta.
\]
Denote the event in the above probability by $\Omega^r_1$. By the
stochastic dynamic equation (\ref{eq:stoc-dym-eqn-B}), we have
\[
\frac{1}{r}\buf(r^2t)(A)=
\frac{1}{r}\sum_{i=B^r(r^2t)+1}^{E^r(r^2 t)}\delta_{v^r_i}(A).
\]
Since $\nu$ is a probability measure on $\R_+$, there exists an
$a>0$ such that \mbox{$\nu(a,\infty)>0$}. We have the following inequality
from the dynamic equation (\ref{eq:mechanism-B}):
%
%
\begin{equation}\label{eq:tech-4-workload}
\frac{1}{r} W^r(r^2t)>a\frac{1}{r}\buf^r(r^2t)(a,\infty)
\ge\frac{a}{r}\sum_{i=B^r(r^2t)+1}^{E^r(r^2 t)}\delta
_{v^r_i}(a,\infty).
\end{equation}
For any $r$, on the event $\Omega^r_1$ there exists a
$t_1\in[0,T+1]$ (random and depending on~$r$) such that
\[
\frac{1}{r}Q^r(r^2t_1)>M.
\]
By (\ref{eq:omegaE-def}), on the event $\Omega^r_E$,
\[
\sup_{t\in[0,T+1]} E^r(r^2t)\le2\lambda r^2(T+1),
\]
for all large enough $r$. Let $M_1=\max(M, 2\lambda T)$ and
$L_1=M$. By (\ref{eq:omega-GC-def}) and (\ref{eq:tech-4-workload}),
on the event
$\Omega^r_{\mathrm{GC}}(M_1,L_1)\cap\Omega^r_E\cap\Omega^r_1$,
\[
\dr W(t_1)>aM\nu^r(a,\infty)-1>aM\nu(a,\infty)-2,
\]
for all large $r$. By (\ref{eq:omegaE}) and (\ref{eq:omega-GC}), we
have for each $M>0$,
\[
\liminf_{r\to\infty}\mathbb{P}^r\Bigl({\sup_{t\in[0,T+1]}\dr
W(t)>aM\nu(a,\infty)-2 }\Bigr)>\eta.
\]
This contradicts the result in Proposition \ref{prop:workload}.
\end{pf}

The following lemma gives a bound on the $(1+p)$th moment of the
measure-valued process, where $p$ is the same as in
conditions (\ref{eq:cond-S-u}) and (\ref{eq:cond-initial-u}).
%
%
\begin{lem}\label{lem:p-moment-bound}
Assume (\ref{eq:cond-A})--(\ref{eq:cond-S-u}),
(\ref{eq:cond-K})--(\ref{eq:cond-initial-u}) and (\ref{eq:cond-HT}).
Fix $T>0$ and $L>1$. For each $\eta>0$, there exists a constant
$M>0$ such that
\[
\liminf_{r\to\infty}\mathbb{P}^r\Bigl({\max_{m\le\lfloor
{rT}\rfloor}\sup
_{t\in[0,L]} \langle{\chi^{1+p}},{\frm\buf(t)+\frm\ser
(t)}\rangle\le M
}\Bigr)\ge1-\eta.
\]
\end{lem}
\begin{pf}
By condition (\ref{eq:cond-initial-u}),
\[
\liminf_{r\to\infty}\mathbb{P}^r\biggl({\biggl\langle{\chi
^{1+p}},{\frac{1}{r}\buf^r(0)+\frac{1}{r}\ser^r(0)}\biggr\rangle\le
\langle{\chi^{1+p}},{\xi^*+\mu^*}\rangle+1}\biggr)=1.
\]
Denote the event in the above by $\Omega^r_0$.
By Lemma \ref{lem:Q-bound}, for any $\eta>0$, there exists a
constant $M'>0$ such that
\[
\liminf_{r\to\infty}\mathbb{P}^r\biggl({\max_{m\le\lfloor
{rT}\rfloor}\sup
_{t\in[0,L]}\frac{1}{r}Q^r(rm+rt)\le M' }\biggr)\ge1-\eta/2.
\]
Denote the event in the above by $\Omega^r_1(M)$. Fix
$M_1=\max(M',\lambda(T+1))$ and $L_1=\lambda(L+1)+2M'$. By
Lemma \ref{lem:glivenko-cantelli},
\[
\lim_{r\to\infty}\mathbb{P}^r({\Omega^r_{\mathrm
{GC}}(M_1,L_1)})=1.
\]
To prove the lemma, it suffices to show that there exists an $M>0$
such that on the event
$\Omega^r_0\cap\Omega^r_1(M')\cap\Omega^r_{\mathrm
{GC}}(M_1,L_1)\cap
\Omega^r_E$,
\[
\max_{m\le\lfloor{rT}\rfloor}\sup_{t\in[0,L]}
\biggl\langle{\chi^{1+p}},{\frac{1}{r}\buf^r(rm+rt)+\frac{1}{r}\ser
^r(rm+rt)}\biggr\rangle
\le M,
\]
for all large $r$. In the remainder of the proof, all random
quantities of the $r$th system is evaluated at a sample path in the
event $\Omega^r_0\cap\Omega^r_1(M')
\cap\Omega^r_{\mathrm{GC}}(M_1,L_1) \cap\Omega^r_E$.

We first find a bound for
$\max_{m\le\lfloor{rT}\rfloor}\sup_{t\in[0,L]}\langle{\chi
^{1+p}},{\frac{1}{r}\buf^r(rm+rt)}\rangle$.
By the dynamic equation (\ref{eq:stoc-dym-eqn-B}), we have that for
all $m\le\lfloor{rT}\rfloor$ and $t\in[0,L]$,
\[
\biggl\langle{\chi^{1+p}},{\frac{1}{r}\buf^r(rm+rt)}\biggr\rangle
=\Biggl\langle{\chi^{1+p}},{\frac{1}{r}\sum
_{B^r(rm+rt)}^{E^r(rm+rt)}\delta_{v^r_i}}\Biggr\rangle.
\]
By (\ref{eq:omegaE-def}) and the definition of $\Omega^r_1(M')$, we
have
%
%
\begin{eqnarray}\label{eq:tech-4-E-bound}
\max_{m\le\lfloor{rT}\rfloor}\sup_{t\in[0,L]} E^r(rm+rt)&<&
\lambda
r^2(T+1)\le
r^2M_1,\\
\label{eq:tech-4-Q-bound}
\max_{m\le\lfloor{rT}\rfloor}\sup_{t\in[0,L]}Q^r(rm+rt)&<& rM'\le rL_1
\end{eqnarray}
for all large enough $r$. So
\[
\max_{m\le\lfloor{rT}\rfloor}\sup_{t\in[0,L]}\biggl\langle{\chi
^{1+p}},{\frac{1}{r}\buf^r(rm+rt)}\biggr\rangle
\le\sup_{-rM_1<n<r^2M_1}\langle{\chi^{1+p}},{\fr\eta
(n,L_1)}\rangle.
\]
By (\ref{eq:set-of-test-fns}) in the remark after
Lemma \ref{lem:glivenko-cantelli}, the function
$\chi^{1+p}\in\testfn$, which appears in the definition of
$\Omega^r_\mathrm{GC}(M_1,L_1)$. So by (\ref{eq:omega-GC-def})
\[
\max_{m\le\lfloor{rT}\rfloor}\sup_{t\in[0,L]}\langle{\chi
^{1+p}},{\frm\buf(t)}\rangle
\le L_1\langle{\chi^{1+p}},{\nu^r}\rangle+1/2.
\]
It then follows from condition (\ref{eq:cond-S-u}) and Theorem 25.12
in \cite{Billingsley1995}, we have
$\langle{\chi^{1+p}},{\nu^r}\rangle\to\langle{\chi^{1+p}},{\nu
}\rangle$ as
$r\to\infty$. Thus,
\[
\max_{m\le\lfloor{rT}\rfloor}\sup_{t\in[0,L]}\langle{\chi
^{1+p}},{\frm\buf(t)}\rangle
\le L_1\langle{\chi^{1+p}},{\nu}\rangle+1,
\]
for all large $r$.

We now look for a bound for
$\max_{m\le\lfloor{rT}\rfloor}\sup_{t\in[0,L]}\langle{\chi
^{1+p}},{\frac{1}{r}\ser^r(rm+rt)}\rangle$.
It follows from the dynamic equation (\ref{eq:mechanism-S}) that
for any $m\le\lfloor{rT}\rfloor$, $t\in[0,L]$ and Borel set
$A\subset\R_+$,
\begin{eqnarray*}
&&\frac{1}{r}\ser^r(rm+rt)(A)\\
&&\qquad=\frac{1}{r}\ser^r(rm+rt)\bigl(A\cap(0,\infty)\bigr)
+\frac{1}{r}\ser^r(rm+rt)(\{0\})\\
&&\qquad=\frac{1}{r}\ser^r(0)\bigl(A\cap(0,\infty)+S^r(0,rm+rt)\bigr)\\
&&\qquad\quad{}+\sum_{j=0}^{m-1}\frac{1}{r}\sum_{i= B^r(r(m-j-1))+1}^{B^r(r(m-j))}
\delta_{v^r_i}\bigl(A\cap(0,\infty)+S^r({\tau^r_i},rm+rt)\bigr)\\
&&\qquad\quad{}+\frac{1}{r}\sum_{i=\bar B^r(rm)+1}^{B^r(rm+rt)}
\delta_{v^r_i}\bigl(A\cap(0,\infty)+S^r({\tau^r_i},rm+rt)\bigr).
\end{eqnarray*}
Given $0\le j\le m-1$, for those $i$'s with $B^r(r(m-j-1))<i\le
B^r(r(m-j))$ we have
\[
{\tau^r_i}\in[r(m-j-1),r(m-j)].
\]
Let's first consider the case where $Z^r(s)>0$ for all
$s\in[0,rm+rt]$. In this case, by (\ref{eq:cumulative}), the
cumulative service amount
\[
S^r({\tau^r_i},rm+rt)\ge
S^r\bigl(r(m-j),rm\bigr)\ge\frac{rj}{K^r}\ge\frac{j}{2K},
\]
for all large $r$ where the last inequality is due to
(\ref{eq:cond-K}). For those $i$'s such that $\tau^r_i$ larger than
$\bar B^r(rm)$, we use the trivial lower bound
$S^r({\tau^r_i},rm+rt)\ge0$. Also take the trivial lower bound
that $S^r(0,rm+rt)\ge0$. Then we have the following inequality on
the $(1+p)$th moment:
%
%
\begin{eqnarray}\label{ineq:1+p-moment-bound-est-1}
&&\biggl\langle{\chi^{1+p}},{\frac{1}{r}\ser^r(rm+rt)}\biggr\rangle\nonumber\\
&&\qquad\le\biggl\langle{\chi^{1+p}},{\frac{1}{r}\ser^r(0)}\biggr\rangle\nonumber\\[-8pt]\\[-8pt]
&&\qquad\quad{}+\sum_{j=0}^{m-1}
\Biggl\langle{\biggl(\biggl(\chi-\frac{j}{2K}\biggr)^+\biggr)^{1+p}},{\frac{1}{r}\sum
_{i= B^r(r(m-j-1))+1}^{B^r(r(m-j))}\delta_{v^r_i}}\Biggr\rangle\nonumber\\
&&\qquad\quad{}+\Biggl\langle{\chi^{1+p}},{\frac{1}{r}\sum
_{i=B^r(rm)+1}^{B^r(rm+rt)}\delta_{v^r_i}}\Biggr\rangle.\nonumber
\end{eqnarray}
Now, we consider the case where there exists an $s\in[0,rm+rt]$ such
that $Z^r(s)=0$. In this case, let $m_0=\min\{0\le j<
m\dvtx\mbox{there exists an } s\in[r(m-j-1),r(m-j)] \mbox{ such that
}Z^r(s)=0\}$. Pick a point $s_0\in[r(m_0-j-1),r(m_0-j)]$ with
$Z^r(s_0)=0$. Then we can replace time 0 in
(\ref{ineq:1+p-moment-bound-est-1}) with $s_0$ and only consider
intervals $[r(m-j-1),r(m-j)]$ with $j\le m_0-1$. So we have
%
%
\begin{eqnarray}\label{ineq:1+p-moment-bound-est-2}
&&\biggl\langle{\chi^{1+p}},{\frac{1}{r}\ser^r(rm+rt)}\biggr\rangle\nonumber\\
&&\qquad\le\biggl\langle{\chi^{1+p}},{\frac{1}{r}\ser^r(s_0)}\biggr\rangle\nonumber\\[-8pt]\\[-8pt]
&&\qquad\quad{}+\sum_{j=0}^{m_0-1}
\Biggl\langle{\biggl(\biggl(\chi-\frac{j}{2K}\biggr)^+\biggr)^{1+p}},{\frac{1}{r}\sum
_{i= B^r(r(m-j-1))+1}^{B^r(r(m-j))}\delta_{v^r_i}}\Biggr\rangle\nonumber\\
&&\qquad\quad{}+\Biggl\langle{\chi^{1+p}},{\frac{1}{r}\sum
_{i=B^r(rm)+1}^{B^r(rm+rt)}\delta_{v^r_i}}\Biggr\rangle.\nonumber
\end{eqnarray}
It is clear that the upper bound in
(\ref{ineq:1+p-moment-bound-est-2}) is less than or equal to the
upper bound in (\ref{ineq:1+p-moment-bound-est-1}). So we only need
to focus on (\ref{ineq:1+p-moment-bound-est-1}) to estimate an upper
bound for the $(1+p)$th moment. By (\ref{eq:tech-4-E-bound}) and
(\ref{eq:tech-4-Q-bound}), for all $m\le\lfloor{rT}\rfloor$, $t\in
[0,L]$ and
all large~$r$,
\begin{eqnarray*}
-rM'&\le& B^r(rj)\le\lambda r^2 (T+1)\le r^2M_1\\
0&\le& B^r(rj+rt)-B^r(rj)\\
&\le&\lambda
r(L+1)+2rM'<rL_1\qquad\mbox{for all }t\in[0,L].
\end{eqnarray*}
It then follows from (\ref{ineq:1+p-moment-bound-est-1}) that
%
%
\begin{eqnarray}\label{eq:tech-4-1+p-bound}
&&\biggl\langle{\chi^{1+p}},{\frac{1}{r}\ser^r(rm+rt)}\biggr\rangle\nonumber\\
&&\qquad\le\biggl\langle{\chi^{1+p}},{\frac{1}{r}\ser^r(0)}\biggr\rangle\nonumber\\[-8pt]\\[-8pt]
&&\qquad\quad{} +\sup_{-rM_1<n<r^2M_1}\sum_{j=0}^{m-1}
\Biggl\langle{\biggl(\biggl(\chi-\frac{j}{2K}\biggr)^+\biggr)^{1+p}},{\fr\eta
(n,L_1)}\Biggr\rangle\nonumber\\
&&\qquad\quad{} +\sup_{-rM_1<n<r^2M_1}\langle{\chi^{1+p}},{\fr\eta
(n,L_1)}\rangle.\nonumber
\end{eqnarray}
The first term on the right-hand side of the above is bounded by
$1+\langle{\chi^{1+p}}$, ${\xi^*+\mu^*}\rangle$ by the definition of
$\Omega^r_0$.
Again, due to (\ref{eq:omega-GC-def}), condition (\ref{eq:cond-S-u})
and Theorem 25.12 in \cite{Billingsley1995}, the third term on the
right-hand side of the above is bounded by
%
%
\begin{equation}
L_1\langle{\chi^{1+p}},{\nu^r}\rangle+1/2\le L_1\langle{\chi
^{1+p}},{\nu}\rangle+1,
\end{equation}
for all large $r$. It now only remains to deal with the second
term on the right-hand side of (\ref{eq:tech-4-1+p-bound}). Let
\[
\bar F^r_n(x)=\fr\eta(n,L_1)((x,\infty)) \qquad\mbox{for all }x\ge0.
\]
The summation in the second term on the right-hand side of
(\ref{eq:tech-4-1+p-bound}) can be upper bounded by
\begin{eqnarray*}
&& \frac{1}{1+p}\sum_{j=1}^{m-1}
\int_{{j}/({2K})}^\infty\biggl(x-\frac{j}{2K}\biggr)^p\bar F^r_n(x)\,dx\\
&&\qquad\le\frac{2K}{1+p}\sum_{j=1}^{m-1}
\int_{({j-1})/({2K})}^{{j}/({2K})}\int_y^\infty(x-y)^p\bar
F^r_n(x)\,dx\,dy\\
&&\qquad\le\frac{1}{1+p}\int_0^\infty\int_y^\infty(x-y)^p\bar
F^r_n(x)\,dx\,dy.
\end{eqnarray*}
Applying Fubini's theorem, the last bound in the above equals
\begin{eqnarray*}
\frac{1}{1+p}\int_0^\infty\int_0^x(x-y)^p\bar F^r_n(x)\,dy\,dx
&=&\frac{1}{(1+p)^2}\int_0^\infty x^{1+p}\bar F^r_n(x)\,dx\\
&=&\frac{2+p}{(1+p)^2}\langle{\chi^{2+p}},{\fr\eta(n,L_1)}\rangle.
\end{eqnarray*}
%
Since the function $\chi^{2+p}\in\testfn$, it again follows from
(\ref{eq:omega-GC-def}), (\ref{eq:cond-S-u}) and Theorem~25.12 in
\cite{Billingsley1995} that the second term on the right-hand side
of (\ref{eq:tech-4-1+p-bound}) is bounded by
\[
\frac{2+p}{(1+p)^2}L_1\langle{\chi^{2+p}},{\nu^r}\rangle
+1/2\le\frac{2+p}{(1+p)^2}L_1\langle{\chi^{2+p}},{\nu}\rangle+1,
\]
for all large $r$. The proof of this lemma is completed by summing
up all these upper bounds.
\end{pf}


The following proposition summarizes the bound estimates in this
section.
\begin{prop}\label{prop:bound}
Assume (\ref{eq:cond-A})--(\ref{eq:cond-S-u}),
(\ref{eq:cond-K})--(\ref{eq:cond-initial-u}) and (\ref{eq:cond-HT}).
For any $\eta>0$, there exists a constant $M>0$ and an event
$\Omega^r_B(M)$ for each index $r$ such that
%
%
\begin{equation}\label{ineq:Omega-B}
\liminf_{r\to\infty}\mathbb{P}^r({\Omega^r_B(M)})\ge
1-\eta,
\end{equation}
and on the event $\omega\in\Omega^r_B(M)$, we have
\begin{eqnarray*}
\max_{m\le\lfloor{rT}\rfloor}\sup_{t\in[0,L]} \frm Q(t)&\le& M,\\
\max_{m\le\lfloor{rT}\rfloor}\sup_{t\in[0,L]} \frm W(t)&\le& M,\\
\max_{m\le\lfloor{rT}\rfloor}\sup_{t\in[0,L]}
\langle{\chi^{1+p}},{\frm\buf(t)+\frm\ser(t)}\rangle&\le& M.
\end{eqnarray*}
\end{prop}
\begin{pf}
The first and the third inequality follow from
Lemmas \ref{lem:Q-bound} and \ref{lem:p-moment-bound}. The second
inequality follows from Proposition \ref{prop:workload}.
\end{pf}

\subsection{Compact containment}\label{sec:compact-containment}

Recall that a set $\K\subset\M$ is relatively compact if
\[
\sup_{\xi\in\K}\xi(\R_+)<\infty,
\]
and there exists a sequence of nested compact sets $J_n\subset\R_+$
such that $\bigcup J_n=\R_+$ and
\[
\lim_{n\to\infty}\sup_{\xi\in\K}\xi(J_n^c)=0,
\]
where $J_n^c$ denotes the complement of $J_n$; see
\cite{Kallenberg1986}, Theorem A7.5. We establish the following
relative compactness property using the bound estimates in
Section \ref{subsec:some-bound-estimates}.
%
%
\begin{lem}\label{lem:compactness}
Assume (\ref{eq:cond-A})--(\ref{eq:cond-S-u}),
(\ref{eq:cond-K})--(\ref{eq:cond-initial-u}) and (\ref{eq:cond-HT}).
Fix $T>0$ and $L>1$. For each $\eta>0$ there exist a constant $M>0$
and a relatively compact set $\K(M)\subset\M$ such that for all
$\omega\in\Omega^r_B(M)$ (which is introduced in
Proposition \ref{prop:bound}) and $r\in\R_+$,
\[
\frm\buf(\omega,t)\in\K(M)\quad\mbox{and}\quad
\frm\ser(\omega,t)\in\K(M)
\]
for all $t\in[0,T]\mbox{ and }m\le\lfloor{rT}\rfloor$.
\end{lem}
\begin{pf}
Let
\[
\mathbf K(M)=\{\xi\in\M\dvtx\xi(\R_+)\le M
\mbox{ and }\xi((n,\infty))\le M/n^{1+p} \}.
\]
Clearly, $\mathbf K(M)$ is a relatively compact set for any constant
$M>0$. Note that $\frm\buf(\omega,t)(\R_+)$ is bounded by $M$ for
all $m\le\lfloor{rT}\rfloor$, $t\in[0,T]$ and $\omega\in\Omega
^r_B(M)$. By the
Markov inequality, for any $t\ge0$, $m\le\lfloor{rT}\rfloor$ and
$\omega\in\Omega^r_B(M)$,
\[
\frm\buf(\omega,t)((n,\infty))
\le\frac{\langle{\chi^{1+p}},{\frm\buf(\omega,t)}\rangle}{n^{1+p}},
\]
which is bounded by $\frac{M}{n^{1+p}}$ by the definition of
$\Omega^r_B(M)$.

Note that $\frm\ser(\omega,t)(\R_+)$ is bounded by $K^r/r$ by the
policy constraint (\ref{eq:policy-S}). By condition
(\ref{eq:cond-K}), $K^r/r\le K+1$ for all large $r$. The same
argument of $\frm\buf(\omega,t)$ applies for $\frm\ser(\omega,t)$.
\end{pf}

\subsection{Asymptotic regularity}\label{subsec:regularity}

A similar result as in this section was proved in \cite{ZDZ2009}.
However, here we consider a much longer time horizon $[0,\lfloor
{rT}\rfloor+L]$
instead of interval $[0,T]$ in \cite{ZDZ2009}. The proof of the
following result use a combination of ideas in \cite{GromollKurk2007}
and \cite{ZDZ2009}.
\begin{lem}\label{lem:asympt-regul}
Assume (\ref{eq:cond-A})--(\ref{eq:cond-HT}). Fix $T>0$ and $L>1$.
For each $\varepsilon,\eta>0$ there exists a $\kappa>0$ (depending on
$\varepsilon$ and $\eta$) such that
%
%
\begin{equation}\label{ineq:asym-reg}
\liminf_{r\to\infty}\mathbb{P}^r\Bigl({\max_{m\le\lfloor
{rT}\rfloor}\sup
_{t\in[0,L]}\sup_{x\in\R_+} \frm\ser(t)([x,x+\kappa])\le
\varepsilon}\Bigr)\ge1-\eta.
\end{equation}
\end{lem}
\begin{pf}
To prove (\ref{ineq:asym-reg}), it suffices to show
\[
\liminf_{r\to\infty}\mathbb{P}^r\Bigl({\sup_{t\in[0,\lfloor
{rT}\rfloor+L]}\sup_{x\in\R_+} \fro{\ser}(t)([x,x+\kappa])\le
\varepsilon
}\Bigr)\ge1-\eta.
\]
First, we have that for any $\varepsilon, \eta>0$, there exists a
$\kappa$ such that
%
%
\begin{equation}\label{ineq:initial-regular}
\liminf_{r\to\infty}\mathbb{P}^r\Bigl({\sup_{x\in\R_+}\fro{\ser}
(0)([x,x+\kappa])\le\varepsilon/2 }\Bigr)\ge1-\eta/2.
\end{equation}
The proof of this inequality, which is based on
(\ref{eq:cond-initial-noatom}), is exactly the same as the proof of
(5.14) in \cite{ZDZ2009}, so we omit it for brevity.

Now, we need to extend this result to the interval $[0,\lfloor
{rT}\rfloor+L]$.
Denote the event in (\ref{ineq:initial-regular}) by $\Omega_1^r$.
Let
%
%
\begin{equation}\label{eq:omega-reg}
\Omega^r_2(M)=\Omega^r_1\cap\Omega^r_E\cap\Omega^r_B(M)
\cap\Omega^r_{\mathrm{GC}}(2\lambda T,2M).
\end{equation}
By (\ref{eq:omegaE}), (\ref{eq:omega-GC}) and (\ref{ineq:Omega-B}),
there exists an $M>0$ such that
\[
\liminf_{r\to\infty}\mathbb{P}^r({\Omega^r_2(M)})\ge
1-\eta.
\]
In the remainder of the proof, all random objects are evaluated at a
fixed sample path in $\Omega^r_2(M)$.

For any $r>0$, $t\in[0,\lfloor{rT}\rfloor+L]$, we define the time
\[
t_0=\sup\bigl\{
\{s\le t\dvtx\langle{1},{\fro{\ser}(s)}\rangle<\varepsilon/4\}\cup\{0\}
\bigr\}
\]
to be the last time before $t$ that the fluid-scaled number of jobs
in service is less than $\varepsilon/4$. Let
\[
t_1=\max\biggl(t_0, t-\frac{4MK}{\varepsilon}\biggr).
\]
We have the following three cases for discussion.

If $t_1=0$, then by (\ref{ineq:initial-regular}) for each $x\in\R_+$
\[
\fro{\ser}(t_1)\bigl([x,x+\kappa]+S^r(rt_1,rt)\bigr)\le\varepsilon/2.
\]
If $t_1=t_0>0$, then for each $\delta\in(0,t_1)$ there exists an
$s\in(t_1-\delta,t_1]$ such that $\fro Z(s)(\R_+)<\varepsilon/4$.
Since we are only concerned with small $\varepsilon$ [which should be
small enough such that $\fro Z(s)<\varepsilon/4<K^r/r$], $\fro Q(s)=0$
by the policy constraint (\ref{eq:policy-S}). Note that
(\ref{eq:B(t)}) implies
%
%
\begin{equation}\label{eq:arrival-E-bound}
\fro B(s,t)\le\fro E(s,t)+\fro Q(s)\qquad\mbox{for all }s\le t.
\end{equation}
By (\ref{eq:omegaE-def}), we have $\fro
B(s,t_1)\le\lambda\delta+\varepsilon/4$ for all large $r$. For any
Borel set $A\subset\R_+$, by the fluid scaled system dynamic
equation (\ref{eq:mechanism-S}),
\begin{eqnarray*}
\fro{\ser}(t_1)(A)&=&\fro{\ser}(t_1)\bigl(A\cap(0,\infty)\bigr)
+\fro{\ser}(t_1)(\{0\})\\
&\le&\fro{\ser}(s)(\R_+)
+\fro B(s,t_1)\le\varepsilon/4+\lambda\delta+\varepsilon/4,
\end{eqnarray*}
which can be made smalled than $\varepsilon/2$ by choosing $\delta$
suitably small.%

If $t_1=t-\frac{2MK}{\varepsilon}>0$, then since the sharing limit is
$K^r$, we have $S^r(rt_1,rt) \ge\frac{4MKr}{\varepsilon
K^r}\ge\frac{2M}{\varepsilon}$ for all large $r$. So
\[
\fro{\ser}(t_1)\bigl([x,x+\kappa]+S^r(rt_1,rt)\bigr)
\le\fro{\ser}(t_1)\biggl(\biggl[\frac{2M}{\varepsilon},\infty\biggr)\biggr)\le\varepsilon/2,
\]
where the last inequality is due to the Markov's inequality and the
definition of $\Omega^r_B(M)$. To summarize, we have
%
%
\begin{equation}
\fro{\ser}(t_1)\bigl([x,x+\kappa]+S^r(rt_1,rt)\bigr)\le\varepsilon/2.
\end{equation}

By the fluid scaled stochastic dynamic equation
(\ref{eq:mechanism-S}),
%
%
\begin{eqnarray}
\label{eq:3}\qquad
\fro{\ser}(t)([x,x+\kappa])
&=&\fro{\ser}(t_1)\bigl([x,x+\kappa]+S^r(rt_1,rt)\bigr)\nonumber\\[-8pt]\\[-8pt]
&&{} +\frac{1}{r}\sum_{i=r\fro B(t_1)+1}^{r\fro B(t)}
\delta_{v^r_i}\bigl([x,x+\kappa]+S^r(\tau_i, rt)\bigr),
\nonumber
\end{eqnarray}
for each $x>0$. When $x=0$, we have
\begin{eqnarray*}
\fro{\ser}(t)([0,\kappa])&=&\fro{\ser}(t)((0,\kappa])+\fro{\ser
}(t)(\{
0\})\\
&=&\fro{\ser}(t_1)\bigl((0,\kappa]+S^r(rt_1,rt)\bigr)\\
&&{} +\frac{1}{r}\sum_{i=r\fro B(t_1)+1}^{r\fro B(t)}
\delta_{v^r_i}\bigl((0,\kappa]+S^r(\tau_i, rt)\bigr),\\
&\le&\fro{\ser}(t_1)\bigl([0,\kappa]+S^r(rt_1,rt)\bigr)\\
&&{} +\frac{1}{r}\sum_{i=r\fro B(t_1)+1}^{r\fro B(t)}
\delta_{v^r_i}\bigl([0,\kappa]+S^r(\tau_i, rt)\bigr).
\end{eqnarray*}
Since all we need is an upper bound estimate, we stick with
(\ref{eq:3}) for analysis. By the choice of $t_1$, the first term
on the right-hand side (\ref{eq:3}) is always upper bounded by
$\varepsilon/2$. Let $I$ denote the second term on the right-hand side
of the proceeding equation. Now it only remains to show that
$I\le\varepsilon/2$.

Let $t_1<t_2<\cdots<t_N=t$ be a partition of the interval
$[t_1,t]$ such that $|t_{j+1}-t_j|<\delta$ for all $j=1,\ldots,
N-1$, where $\delta$ is to be chosen below.
Note that by the definition of $t_1$,
%
%
\begin{equation}\label{eq:tech-4-N-bound}
N\le\frac{4MK}{\delta\varepsilon}.
\end{equation}
Write $I$ as the summation
\[
I=\sum_{j=1}^{N-1}\frac{1}{r}\sum_{i=r\fro B(t_j)+1}^{r\fro
B(t_{j+1})} \delta_{v^r_i}\bigl([x,x+\kappa]+S^r(\tau_i, rt)\bigr).
\]
Recall that $\tau^r_i$ is the time that the $i$th job starts
service, so $rt_j\le\tau^r_i\le rt_{j+1}$ for those $i\in[r\fro
B(t_j)+1,r\fro B(t_{j+1})]$. This implies that
\[
S^r(rt_{j+1},rt)\le S^r(\tau_i, rt)\le S^r(rt_j,rt).
\]
By the definition of $t_1$, we have $\fro Z(s)\ge\varepsilon/4$ for all
$s\in[t_1,t]$. This gives
\[
S^r(rt_j,rt_{j+1})\le\frac{4\delta}{\varepsilon}.
\]
So for any $i\in[r\fro B(t_j)+1,r\fro B(t_{j+1})]$, we have
$[x,x+\kappa]+S^r(\tau_i, rt)\subseteq C_j$, where
\[
C_j=\biggl[x+S^r(rt_{j+1},rt),
x+\kappa+S^r(rt_{j+1},rt)+\frac{4\delta}{\varepsilon}\biggr].
\]
Thus,
\[
I\le\sum_{j=1}^{N-1}\frac{1}{r}
\sum_{i=r\fro B(t_j)+1}^{r\fro B(t_{j+1})}\delta_{v^r_i}(C_j).
\]
By (\ref{eq:omegaE-def}), (\ref{eq:arrival-E-bound}) and the
definition of $\Omega^r_B(M)$, we have for all $j=1,\ldots, N-1$
\begin{eqnarray*}
-rM&\le& r\fro B(t_j)\le r(\lambda rT+L+1+M)\le r^22\lambda T,\\
\fro B(t_j,t_{j+1})&\le&\lambda\delta+1+M\le2M,
\end{eqnarray*}
where the last inequality in each of the above bound holds because
we only care about small $\delta$ and large $r$. Choose
$\varepsilon_1=\frac{\delta\varepsilon^2}{16MK}$. By
(\ref{eq:omega-GC-def}),
\[
\Biggl|
\frac{1}{r}\sum_{i=r\fro B(t_j)+1}^{r\fro B(t_{j+1})}
\delta_{v^r_i}(C_j)-\bigl(\fro B(t_{j+1})-\fro B(t_j)\bigr)\nu^r(C_j)
\Biggr|\le\varepsilon_1,
\]
for all large $r$. This implies that
\[
I\le
\sum_{j=1}^{N-1}[\fro B(t_{j+1})-\fro B(t_j)]\nu^r(C_j)+N\varepsilon_1.
\]
Let
$\varepsilon_2=\frac{\varepsilon}{8(\lambda{4MK}/{\varepsilon}+M+1)}$.
Since $C_j^{\kappa\varepsilon_2}$ is a close interval with length
$\kappa+\frac{4\delta}{\varepsilon}+\kappa\varepsilon_2$, by condition
(\ref{eq:cond-nu-noatom}), we can choose $\kappa,\delta<1$ small
enough such that
\[
\nu(C_j^{\kappa\varepsilon_2})\le\varepsilon_2,
\]
where $C_j^{\varepsilon_2}$ is the $\varepsilon_2$-enlargement of $C_j$.
By (\ref{eq:cond-S}), we also have
\[
\nu^r(C_j)\le\nu(C_j^{\kappa\varepsilon_2})+\kappa\varepsilon_2
\le\nu(C_j^{\kappa\varepsilon_2})+\varepsilon_2,
\]
for all large enough $r$. Thus, we conclude that
\begin{eqnarray*}
I&\le&2\varepsilon_2\sum_{j=1}^{N-1}[\fro B(t_{j+1})-\fro
B(t_j)]+N\varepsilon
_1\\
&\le&2\varepsilon_2[\fro B(t)-\fro B(t_1)]+N\varepsilon_1\\
&\le&2\varepsilon_2\biggl(\lambda\frac{4MK}{\varepsilon}+1+M\biggr)+N\varepsilon_1,
\end{eqnarray*}
where the last inequality is again due to (\ref{eq:omega-reg}),
(\ref{eq:arrival-E-bound}) and (\ref{eq:tech-4-N-bound}). Finally,
by the choice of $\varepsilon_1, \varepsilon_2$, we have that
$I\le\varepsilon/2$.
\end{pf}

\subsection{Oscillation bound}\label{sec:oscillation-bound}

Consider a c\`adl\`ag function $\zeta(\cdot)$ on a fixed interval $[0,L]$
taking values in a metric space $(\mathbf E,\gd)$. The \textit{modulus
of continuity} is defined to be
\[
\osc{\zeta}{\delta}=\sup_{s,t\in[0,L],|s-t|<\delta}\gd[\zeta
(s),\zeta(t)].
\]
If the metric space is $\R$, we just use the Euclidian norm; if the
space is $\M$ or $\M\times\M$, we use the Prohorov metric $\pov$
defined in Section \ref{sec:introduction}. We have the following
bound on the oscillation of the shifted fluid scaled measure-valued
processes.
\begin{lem}\label{lem:oscillation-bound}
Assume (\ref{eq:cond-A})--(\ref{eq:cond-HT}). Fix $T>0$ and $L>1$.
For each $\varepsilon,\eta>0$ there exists a $\delta>0$ (depending on
$\varepsilon$ and $\eta$) such that
%
%
\begin{equation}\label{eq:oscillation-bound}
\liminf_{r\to\infty}\mathbb{P}^r\Bigl({\max_{m\le\lfloor
{rT}\rfloor}\max
( \osc{\frm\buf}{\delta},\osc{\frm\ser}{\delta} ) \le
\varepsilon}\Bigr)\ge1-\eta.\hspace*{-28pt}
\end{equation}
\end{lem}

The proof of this lemma, which builds on the asymptotic regularity in
Lem\-ma~\ref{lem:asympt-regul}, uses the exactly same argument as in the
proof of Lemma 5.6 based on Lemma 5.5 in \cite{ZDZ2009}. We omit
this proof for brevity.

Fix $T>0$ and $L>0$. For any sequence 
$\{\delta_i\}$, consider the following set:
%
%
\begin{equation}\label{eq:omega-reg-osc}
\biggl\{
\max_{m\le\lfloor{rT}\rfloor}\max(
\osc{\fr\buf}{\delta_j},\osc{\fr\ser}{\delta_j}
)\le\frac{1}{j}
\biggr\}.
\end{equation}
Denote the 
sequence 
$\{\delta_i\}$ by $\mathcal S$. To emphasize the dependency on
$\mathcal S$ and $j$, denote the above event by $\Omega^r_R(\mathcal
S,j)$. By Lemmas \ref{lem:asympt-regul} and
\ref{lem:oscillation-bound}, for any $\eta>0$, there exists an
$\mathcal S$ such that
%
%
\begin{equation}\label{ineq:omega-reg-osc}
\liminf_{r\to\infty}\mathbb{P}^r({\Omega^r_R(\mathcal
S,j)})
\ge1-\frac{\eta/2}{2^j}\qquad
\mbox{for }j=1,2,\ldots.
\end{equation}
%
This implies that for any finite number $n\in\N$, we have
\[
\liminf_{r\to\infty}\mathbb{P}^r\Biggl({\bigcap_{j=1}^{n}\Omega
^r_R(\mathcal S,j)}\Biggr)
\ge1-{\eta/2}.
\]
Let $r(n)$ denote the smallest number such that
%
%
\begin{equation}\label{ineq:omega-S-j-r-j}
\mathbb{P}^r\Biggl({\bigcap_{j=1}^{n}\Omega^r_R(\mathcal S,j)}
\Biggr)\ge1-\eta\qquad
\mbox{for all }r\ge r(n).
\end{equation}
For any $n\le n'$, we have $r(n)\le r(n')$ since
$\bigcap_{j=1}^{n}\Omega^r_R(\mathcal S,j) \supseteq
\bigcap_{j=1}^{n'}\Omega^r_R(\mathcal S,j)$. Let
\[
n(r)=\sup\bigl\{\{n\in\Z_+\dvtx r(n)\le r\}\cup\{0\}\bigr\}.
\]
[By this definition, we allow $n(\cdot)$ to be infinite. For example,
when the function $r(\cdot)$ has an upper bound. In fact, $n(\cdot)$
can be viewed as the ``inverse'' of $r(\cdot)$.] It is clear that
$n(\cdot)$ is nondecreasing. Note that for any $n_0>0$ there exists
$r_0=r(n_0)$ such that $n(r)\ge n_0$ for all $r\ge r_0$. Thus, we have
that
\[
\lim_{r\to\infty}n(r)=\infty.
\]
Now define
%
%
\begin{equation}\label{eq:omegaR}
\Omega^r_R(\mathcal S)=\bigcap_{j=1}^{n(r)}\Omega^r_R(\mathcal S,j).
\end{equation}
Note that $\Omega^r_R(\mathcal S)$ is not empty for all large enough
$r$ [since $n(r)>1$ for all large enough $r$], and in this case,
$\mathbb{P}^r({\Omega^r_R(\mathcal S)})\ge1-\eta$.
So we conclude that
%
%
\begin{equation}\label{ineq:omegaR}
\liminf_{r\to\infty}\mathbb{P}^r({\Omega^r_R(\mathcal S)}
)\ge1-\eta.
\end{equation}
%
Denote
%
%
\begin{equation}\label{eq:Omega-M}
\Omega^r(M,\mathcal S)
=
\Omega^r_B(M)
\cap\Omega^r_R(\mathcal S).
\end{equation}
%
For any $r$, the $r$th system is defined on the probability space
$(\Omega^r,\mathbb P^r,\mathcal F^r)$. The stochastic processes
$\buf^r(\cdot)$ and $\ser^r(\cdot)$ are actually measurable functions
on $\Omega^r$. From now on, we explicitly write these processes down
in the form of $\buf^r(\omega,\cdot)$ and $\ser^r(\omega,\cdot)$ to
indicate that they are evaluated on the sample path
$\omega\in\Omega^r$. We are now ready to present the precompactness
result.
\begin{theorem}\label{thm:precompactness}
Assume (\ref{eq:cond-A})--(\ref{eq:cond-HT}). Fix $T>0$ and $L>1$.
For each $\eta>0$, the exists a constant $M>0$ and an $\mathcal S$
such that such that
%
%
\begin{equation}\label{ineq:Omega-M}
\liminf_{r\to\infty}\mathbb{P}^r({\Omega^r(M,\mathcal S)}
)\ge1-\eta.
\end{equation}
Suppose $\{r_n\}_{n\in\N}$ is a sequence in $\R_+$ which goes to
infinity. Any sequence of functions
$\{(\bar\buf^{r_n, m_n}(\omega_n,\cdot),
\bar\ser^{r_n,m_n}(\omega_n,\cdot))\}_{n\in\N}$
with $\omega_n\in\Omega^{r_n}(M,\mathcal S)$ and
$m_n\le\lfloor{r_nT}\rfloor$ for each $n\in\N$ has a subsequence
$\{(\bar\buf^{r_{n_i},m_{n_i}}(\omega_{n_i},\cdot),
\bar\ser^{r_{n_i},m_{n_i}}(\omega_{n_i},\cdot))\}_{i\in\N}$
such that
\[
\upsilon_L[
(\bar\buf^{r_{n_i},m_{n_i}}(\omega_{n_i},\cdot),
\bar\ser^{r_{n_i},m_{n_i}}(\omega_{n_i},\cdot)),
(\tilde\buf(\cdot),\tilde\ser(\cdot))
]\to0\qquad \mbox{as }i\to\infty,
\]
for some process $(\tilde\buf(\cdot),\tilde\ser(\cdot))$ which is
continuous, where $\upsilon_L$ is the uniform metric defined in
(\ref{eq:sup-norm}).
\end{theorem}
\begin{pf}
For a fixed $\eta>0$, pick an $M>0$ that satisfies
(\ref{ineq:Omega-B}) and construct an $\mathcal S$ so that it
satisfies (\ref{ineq:omega-reg-osc}). Define $\Omega^r(M,\mathcal
S)$ via (\ref{eq:Omega-M}). The probability inequality
(\ref{ineq:Omega-M}) follows immediately from (\ref{ineq:Omega-B})
and (\ref{ineq:omegaR}). The space $\M\times\M$ endowed with the
metric $\pov$ (defined in Section \ref{subsec:notatioin}) is
complete. Lemma \ref{lem:compactness} verifies condition ($a$) in
Theorem 3.6.3 of \cite{EthierKurtz1986}. For any $\varepsilon>0$
there exists a $j_0$ such that $1/j<\varepsilon$ for all $j>j_0$. By
(\ref{eq:omega-reg-osc}) and (\ref{eq:omegaR}), we have that when
$\delta\le\delta_{j_0}$ and $r>r(j_0)$, where $\delta_{j_0}$ is
specified in $\mathcal S$ and $r(n)$ is defined in
(\ref{ineq:omega-S-j-r-j}),
%
%
\begin{equation}\label{eq:tech-5-osc}
\max(\oscw{\frm\buf(\omega^r,\cdot)}{\delta},
\oscw{\frm\ser(\omega^r,\cdot)}{\delta})\le\varepsilon,
\end{equation}
for any $\omega^r\in\Omega^r(M,\mathcal S)$ and $m\le\lfloor
{rT}\rfloor$. This
verifies condition ($b$) in Theorem 3.6.3 of
\cite{EthierKurtz1986}. So the sequence
$\{(\bar\buf^{r_n,m_n}(\omega^{r_n},\cdot),
\bar\ser^{r_n,m_n}(\omega^{r_n},\cdot))\}_{n\in\N}$ is precompact in
the space $\D([0,T], \M\times\M)$ endowed with the Skorohod $J_1$
topology. In other words, there is a convergent subsequence. The
limit of this subsequence is continuous by the oscillation bound
(\ref{eq:tech-5-osc}). So convergence in the Skorohod
$J_1$-topology is the same as convergence in the uniform metric
defined in Section \ref{subsec:notatioin}.
\end{pf}

\section{State-space collapse}\label{sec:diff-appr}

In this section, we establish the state-space collapse
(Theorem \ref{thm:ssc}). The task is divided into the following
steps: we first show that the limits in
Theorem \ref{thm:precompactness}, which called fluid limits, are fluid
model solutions; the set of fluid limits is ``rich'' in the sense that
itself and the set of shifted fluid scaled process mutually
approximates each other (Lemmas \ref{lem:asymp-close} and
\ref{lem:fluid-unif-approx}); the proof of the state-space collapse
result is finally presented based on the richness of fluid limits and
the properties of fluid model solution
(Theorems \ref{thm:convergence-to-invariant-mani} and
\ref{thm:invariant-manifold}).

\subsection{Fluid limits}\label{subsec:fluid-limits}
Let $\ds_L(M,\mathcal S)$ denote the set of fluid limits of all
convergent subsequences of sequences in
Theorem \ref{thm:precompactness}. It is then quite clear that we have
the following property.
\begin{lem}\label{lem:asymp-close}
Assume (\ref{eq:cond-A})--(\ref{eq:cond-HT}). The set of fluid
limits $\ds_L(M,\mathcal S)$ is nonempty. Pick an element
$(\tilde\buf(\cdot),\tilde\ser(\cdot))\in\ds_L(M,\mathcal S)$, for
any $\varepsilon>0$ and $r_0\in\R_+$, there exists an $r\ge r_0$,
$m\le\lfloor{rT}\rfloor$ and $\omega\in\Omega^r(M,\mathcal S)$
such that
\[
\upsilon_L[
(\frm\buf(\omega,\cdot),\frm\ser(\omega,\cdot)),
(\tilde\buf(\cdot),\tilde\ser(\cdot))
]\le\varepsilon.
\]
\end{lem}

Roughly speaking, this lemma says that any element in
$\ds_L(M,\mathcal S)$ can be approximated by a shifted fluid scaled
process of the $r$th system evaluated at some sample path in
$\Omega^r(M,\mathcal S)$ with arbitrarily large index $r$. This helps
prove the following property of the fluid limits.

Fix a constant $0<q<p$, where $p$ is the same one as in
(\ref{eq:cond-S-u}) and (\ref{eq:cond-initial-u}). Recall the subset
$\init^q_{3M}$ of all valid initial conditions defined in
(\ref{eq:comp-set-init}).
\begin{lem}\label{lem:prop-fluid-limit}
Assume (\ref{eq:cond-A})--(\ref{eq:cond-HT}). Fix $L>0$ and
$0<q<p$. Any element
$(\tilde\buf(\cdot),\tilde\ser(\cdot))\in\ds_L(M,\mathcal S)$
is a
critically loaded fluid model solution with initial condition
belongs to $\init^q_{3M}$.
\end{lem}
\begin{pf}
We first show that the initial condition
$(\tilde\buf(0),\tilde\ser(0))\in\init^q_{3M}$. By the definition
of the fluid limit, there exists a subsequence
\[
(\bar\buf^{r_i,m_i}(\omega_i,0),\bar\ser^{r_i,m_i}(\omega_i,0))
\to(\tilde\buf(0),\tilde\ser(0))\qquad \mbox{as }i\to\infty,
\]
where the above convergence is in the Prohorov metric. By
Proposition \ref{prop:bound} and the LPS policy, we have
\begin{eqnarray*}
\langle{1},{\bar\buf^{r_i,m_i}(\omega_i,0)+\bar\ser
^{r_i,m_i}(\omega_i,0)}\rangle
&<&M+K+1,\\
\langle{\chi^{1+p}},{\bar\buf^{r_i,m_i}(\omega_i,0), \bar\ser
^{r_i,m_i}(\omega_i,0)}\rangle
&<&M,
\end{eqnarray*}
for all large $r_i$. This implies that for any $0\le q<p$,
\begin{eqnarray*}
&&\hspace*{-5pt} \langle{\chi^{1+q}},{\bar\buf^{r_i,m_i}(\omega_i,0) +\bar
\ser^{r_i,m_i}(\omega_i,0)}\rangle\\
&&\hspace*{-5pt}\qquad\le\langle{1},{\bar\buf^{r_i,m_i}(\omega_i,0) +\bar\ser
^{r_i,m_i}(\omega_i,0)}\rangle
+\langle{\chi^{1+p}},{\bar\buf^{r_i,m_i}(\omega_i,0) +\bar\ser
^{r_i,m_i}(\omega_i,0)}\rangle\\
&&\hspace*{-5pt}\qquad< 2M+K+1.
\end{eqnarray*}
By the corollary of Theorem 25.12 in \cite{Billingsley1995}, we
have that for any $0\le q<p$,
\[
\langle{\chi^{1+q}},{\bar\buf^{r_i,m_i}(\omega_i,0) +\bar\ser
^{r_i,m_i}(\omega_i,0)}\rangle
\to\langle{\chi^{1+q}},{\tilde\buf(0)+\tilde\ser(0)}\rangle\qquad
\mbox{as } i\to\infty.
\]
Since we can take $M$ big enough such that $M>K+1$, this implies that
$\langle{\chi^{1+q}},{\tilde\buf(0)+\tilde\ser(0)}\rangle<3M$ and
$\langle{\chi},{\tilde\buf(0)+\tilde\ser(0)}\rangle<3M$, which
yields the
result.

By Lemma \ref{lem:asymp-close}, any fluid limit
$(\tilde\buf(\cdot),\tilde\ser(\cdot))$ can be approximated by a
shifted fluid scaled process of the $r$th system evaluated at some
sample path in $\Omega^r(M,\mathcal S)$ with arbitrarily large index
$r\in\R_+$; the state descriptor of the $r$th system satisfies the
stochastic dynamic equations (\ref{eq:stoc-dym-eqn-B}) and
(\ref{eq:stoc-dym-eqn-S}). It then follows from the same argument
as in Lemmas 6.1 and 6.2 in \cite{ZDZ2009} that each fluid limit
satisfies the fluid model equations (\ref{eq:fluid-dym-eqn-B}) and
(\ref{eq:fluid-dym-eqn-S}) and constraints
(\ref{constr:non-decr})--(\ref{constr:policy-S}). In fact,
\cite{ZDZ2009} is more general in the sense that the traffic
intensity is allowed to be any positive number instead of being 1 as
required in this paper.
\end{pf}

\subsection{Uniform approximation}\label{subsec:unif-approx}

Lemma \ref{lem:fluid-unif-approx} in the following is analogous to
Lemma 4.1 in \cite{Bramson1998}. In contrast to
Lemma \ref{lem:asymp-close} above, this lemma says that any shifted
fluid scaled process of the $r$th system evaluated at some sample path
in $\Omega^r(M,\mathcal S)$ with index $r$ large enough can be
approximated by some element in $\ds_L(M,\mathcal S)$, which has been
proved to be a fluid model solution in
Lemma \ref{lem:prop-fluid-limit}. This result will help prove the
state-space collapse result for diffusion scaled processes.
\begin{lem}\label{lem:fluid-unif-approx}
Assume (\ref{eq:cond-A})--(\ref{eq:cond-HT}). For each
$\varepsilon>0$, there exists an $r_0\in\R_+$ such that for any $r\ge
r_0$, $m\le\lfloor{rT}\rfloor$ and $\omega\in\Omega^r(M,\mathcal
S)$, we can\vspace*{1pt}
find a $(\tilde\buf(\cdot),\tilde\ser(\cdot))\in\ds
_L(M,\mathcal S)$
satisfying
\[
\upsilon_L[(\frm\buf(\omega,\cdot),\frm\ser(\omega,\cdot))
,(\tilde\buf(\cdot),\tilde\ser(\cdot))]<\varepsilon.
\]
\end{lem}
\begin{pf}
Assume it is not true. Then there exists an $\varepsilon>0$ such that
for any natural number $i$ there exist an $r_i>i$, $m_i\in\lfloor
{rT}\rfloor$
and $\omega_i\in\Omega^r(M,\mathcal S)$ such that
\[
\upsilon_L[(\bar\buf^{r_i,m_i}(\omega_i,\cdot),
\bar\ser^{r_i,m_i}(\omega_i,\cdot))
,(\tilde\buf(\cdot),\tilde\ser(\cdot))]\ge\varepsilon,
\]
for all $(\tilde\buf(\cdot),\tilde\ser(\cdot))\in\ds
_L(M,\mathcal
S)$. However, by Theorem \ref{thm:precompactness}, the sequence
\[
\{(\bar\buf^{r_i,m_i}(\omega_i,\cdot),
\bar\ser^{r_i,m_i}(\omega_i,\cdot))\}_{i=0}^\infty
\]
contains a convergent subsequence, the limit of which must be in
$\ds_L(M,\mathcal S)$. This is a contradiction.
\end{pf}

\subsection{Proof of state-space collapse}
\label{subsec:proof-ssc}

With all the preparation, we finally present the proof of state-space
collapse.
\begin{pf*}{Proof of Theorem \ref{thm:ssc}}
By (\ref{ineq:Omega-M}), it suffices to show that for each
$\varepsilon>0$, there exists an $r_0$ such that when $r>r_0$,
%
%
\begin{equation}\label{eq:ssc-suff-1}
\sup_{\omega\in\Omega^r(M,\mathcal S)}\sup_{t\in[0,T]}
\pov[(\dr\buf(\omega,t),\dr\ser(\omega,t)),
\Delta_{K,\lambda}\dr W(\omega,t)]<\varepsilon.
\end{equation}
In the following, we fix $r>r_0$ and $\omega\in\Omega^r(M,\mathcal
S)$. By Lemma \ref{lem:prop-fluid-limit}, any
$(\tilde\buf(\cdot)$, $\tilde\ser(\cdot))\in\ds_L(M,\mathcal S)$
is a
critically loaded fluid model solution with initial condition
$(\xi,\mu)\in\init^q_{3M}$. Denote
\[
\tilde W(\cdot)=\langle{\chi},{\tilde\buf(\cdot)+\tilde\ser
(\cdot)}\rangle.
\]
It follows from the workload conservation property
(\ref{eq:prop-workload}) that $\tilde W(\cdot)\equiv
\langle{\chi},{\xi+\mu}\rangle$. By
Theorem \ref{thm:convergence-to-invariant-mani}, there exists an
$L^*>0$ such that when $s>L^*$,
%
%
\begin{equation}\label{eq:tech-6-invar}
\pov[(\tilde\buf(s),\tilde\ser(s)),\Delta_{K,\nu} \tilde
W(s)]<\varepsilon/3,
\end{equation}
for all $(\tilde\buf(\cdot),\tilde\ser(\cdot))\in\ds
_L(M,\mathcal
S)$. Now, fix a constant $L>L^*+1$. Note that
\[
[0,r^2T]\subset[0,rL^*]\cup\bigcup_{m=0}^{\lfloor{rT}\rfloor
}[r(m+L^*),r(m+L)].
\]
By the definition of diffusion and shifted fluid scaling, to show
(\ref{eq:ssc-suff-1}) it suffices to show
%
%
\begin{eqnarray}
\label{eq:ssc-suff-2-m}
\max_{m\le\lfloor{rT}\rfloor}\sup_{s\in[L^*,L]}
\pov[(\frm\buf(\omega,s),\frm\ser(\omega,s)),
\Delta_{K,\lambda}\frm W(\omega,s)]&<&\varepsilon,\\
\label{eq:ssc-suff-2-0}
\sup_{s\in[0,L^*]}
\pov[(\fro\buf(\omega,s),\fro{\ser}(\omega,s)),
\Delta_{K,\lambda}\fro W(\omega,s)]&<&\varepsilon.
\end{eqnarray}

We first prove (\ref{eq:ssc-suff-2-m}). Fix an $m\le\lfloor
{rT}\rfloor$. By
Lemma \ref{lem:fluid-unif-approx}, for any $\varepsilon'>0$, there
exists a $(\tilde\buf(\cdot),\tilde\ser(\cdot))\in\ds
_L(M,\mathcal
S)$ (depending on $r,m$ and $\omega$) such that
%
%
\begin{equation}\label{eq:tech-6-immed-01}
\upsilon_L[(\frm\buf(\omega,\cdot),\frm\ser(\omega,\cdot)),
(\tilde\buf(\cdot),\tilde\ser(\cdot))]<\varepsilon'.
\end{equation}
By the definition of $\Omega^r(M,\mathcal S)$ and
Proposition \ref{prop:bound}, following the same proof as in
Lemma \ref{lem:prop-fluid-limit}, we have that for each fixed
$0<q<p$,
\begin{eqnarray*}
\langle{\chi^{1+q}},{\tilde\buf(t)+\tilde\ser(t)}\rangle&<&3M,\\
\langle{\chi^{1+q}},{\frm\buf(\omega,t)+\frm\ser(\omega
,t)}\rangle&<&3M,
\end{eqnarray*}
for all $t\in[0,L]$. It then follows from
Lemma \ref{lem:workload-pov} and by taking $\varepsilon'$ small enough
that
%
%
\begin{equation}\label{eq:tech-6-immd-2}
{\sup_{t\in[0,L]}}|\tilde W(t)-\frm W(\omega,t)|
<\frac{\varepsilon}{3 \max({1}/{\beta},{1}/{\beta_e})}.
\end{equation}
Note that for any real numbers $w_1,w_2$, by the definition of the
lifting map $\Delta_{K,\nu}$ and the metric $\pov$, we have
\begin{eqnarray*}
&& \pov[\Delta_{K,\nu}w_1,\Delta_{K,\nu}w_2]\\
&&\qquad<\max\biggl(\pov\biggl[\frac{(w_1-K\beta_e)^+}{\beta}\nu,
\frac{(w_2-K\beta_e)^+}{\beta}\nu\biggr],\\
&&\qquad\hspace*{64.35pt}\pov\biggl[\frac{w_1\wedge K\beta_e}{\beta_e}\nu_e,
\frac{w_2\wedge K\beta_e}{\beta_e}\nu_e\biggr]
\biggr).
\end{eqnarray*}
It is clear that for any $a,b\ge0$ and Borel set $A\subseteq\R$,
we have that $a\nu(A)\le b\nu(A)+|b-a|\le b\nu(A^{|b-a|})+|b-a|$,
where $A^{|b-a|}$ is the $|b-a|$-enlargement of $A$. Similarly, we
have $b\nu(A)\le b\nu(A^{|b-a|})+|b-a|$. So $\pov[a\nu,b\nu]\le
|b-a|$. This implies that
\begin{eqnarray*}
\pov\biggl[\frac{(w_1-K\beta_e)^+}{\beta}\nu,
\frac{(w_2-K\beta_e)^+}{\beta}\nu\biggr]
&\le&\biggl|\frac{(w_1-K\beta_e)^+}{\beta}
-\frac{(w_2-K\beta_e)^+}{\beta}\biggr|
\\
&\le&\frac{1}{\beta}|w_1-w_2|.
\end{eqnarray*}
Following the same argument, we have
\[
\pov\biggl[\frac{w_1\wedge K\beta_e}{\beta_e}\nu_e,
\frac{w_2\wedge K\beta_e}{\beta_e}\nu_e\biggr]
\le\biggl|\frac{w_1\wedge K\beta_e}{\beta_e}
-\frac{w_2\wedge K\beta_e}{\beta_e}\biggr|
\le\frac{1}{\beta_e}|w_1-w_2|.
\]
Thus, we conclude that
%
%
\begin{equation}\label{eq:tech-6-immd-3}
\pov[\Delta_{K,\nu}w_1,\Delta_{K,\nu}w_2]
<\max\biggl(\frac{1}{\beta},\frac{1}{\beta_e}\biggr)|w_1-w_2|.
\end{equation}
So (\ref{eq:ssc-suff-2-m}) follows from (\ref{eq:tech-6-invar}) and
(\ref{eq:tech-6-immed-01})--(\ref{eq:tech-6-immd-3}).

It now remains to show (\ref{eq:ssc-suff-2-0}). By
Lemma \ref{lem:fluid-unif-approx}, for any $\varepsilon'>0$, there
exists a $(\tilde\buf(\cdot),\tilde\ser(\cdot))\in\ds
_L(M,\mathcal
S)$ (depending on $r$ and $\omega$) such that
%
%
\begin{equation}\label{eq:tech-6-immed-1}
\upsilon_L[(\fro\buf(\omega,\cdot),\fro{\ser}(\omega,\cdot)),
(\tilde\buf(\cdot),\tilde\ser(\cdot))]<\varepsilon'.
\end{equation}
By conditions (\ref{eq:cond-initial}) and (\ref{eq:cond-init-ssc}),
we have that
\[
(\tilde\buf(0),\tilde\ser(0))=\Delta_{K,\nu}\tilde W(0).
\]
In
other words, the initial condition $(\tilde\buf(0),\tilde\ser(0))$
is an equilibrium state. Since
$(\tilde\buf(\cdot),\tilde\ser(\cdot))$ is a fluid model solution,
by Theorem \ref{thm:invariant-manifold},
\[
(\tilde\buf(t),\tilde\ser(t))=\Delta_{K,\nu}\tilde W(t)\qquad
\mbox{for all }t\ge0.
\]
So (\ref{eq:ssc-suff-2-0}) follows immediately from
(\ref{eq:tech-6-immd-2})--(\ref{eq:tech-6-immed-1}) and the above
equation.\vadjust{\goodbreak}
\end{pf*}

\begin{appendix}

\section{An integration by parts formula for Lebesgue--Stieltjes
integral}
\label{sec:int-parts-LS}

The following lemma is used in the derivation of (\ref{eq:key-pre})
and in the proof of Lemma \ref{lem:Ax-converges}. We do not require
the continuity of distribution function $F$.
\begin{lem}\label{lem:LS-int-by-parts}
Suppose that $F$ is a probability distribution function with
$F(0)=0$, and $q\in\D([0, \infty), \R)$ has bounded variation on
$[0,b]$ for each $b>0$. 
For any $u>0$,
\[
\int_{[0,u]}[1-F(u-v)]\,dq(v)
=q(u)-[1-F(u)]q(0)-\int_{[0,u]}q(u-v)\,dF(v).
\]
\end{lem}
\begin{pf}
Let $u>0$ be fixed and let $I=[0, u]$. Define $f(v)=1-F(u-v)$. Then
$f$ is a left continuous function in on $(-\infty,u]$. Clearly,
both $f$ and $q$ are functions with bounded variation on the
interval $I$. Let $S$ denote the set of points in $I$ where both $f$
and $q$ are discontinuous. According to Theorem 6.2.2 in
\cite{CarterVanBrunt2000},
%
%
\begin{eqnarray}
\label{eq:int-by-parts}
\int_If\,dq+\int_Iq\,df&=&f(u^+)q(u^+)-f(0^-)q(0^-)+\sum_{a\in
S}A(a)\nonumber\\[-8pt]\\[-8pt]
&=& q(u)-[1-F(u)]q(0)+\sum_{a\in S}A(a),
\nonumber
\end{eqnarray}
where
\begin{eqnarray*}
A(a)&=&\bigl[f(a)-\tfrac{1}{2}\bigl(f(a^+)+f(a^-)\bigr)\bigr]\bigl(q(a^+)-q(a^-)\bigr)\\
&&{} +\bigl[q(a)-\tfrac{1}{2}\bigl(q(a^+)+q(a^-)\bigr)\bigr]\bigl(f(a^+)-f(a^-)\bigr).
\end{eqnarray*}
Since $f$ is continuous on the left and $q$ is continuous on the
right at all $a\in S$, then
%
%
\begin{eqnarray}
\label{eq:int-by-part-adj}
A(a)&=&\bigl[-\tfrac{1}{2}\bigl(f(a^+)-f(a)\bigr)\bigr][q(a)-q(a^-)]\nonumber\\
&&{}+ \bigl[\tfrac{1}{2}\bigl(q(a)-q(a^-)\bigr)\bigr][f(a^+)-f(a)]\\
&=&0.\nonumber
\end{eqnarray}
Now the lemma follows from (\ref{eq:int-by-parts}),
(\ref{eq:int-by-part-adj}) and
\[
\int_Iq\,df = \int_{[0, u]} q(v) \,dF(u-v) = \int_{[0, u]} q(u-v)\,dF(v).
\]
\upqed\end{pf}

\section{A key renewal theorem with uniform convergence}
\label{sec:key-renewal-theorem}

The following result is similar as the key renewal theorem. But the
convergence is shown to be uniform on a set of functions $\mathscr
H$ as specified below.
\begin{lem}\label{lem:key-renewal-thm-unif}
Assume that each $h\in\mathscr H$ is nonnegative and
nonincreasing. Assume that
%
%
\begin{eqnarray}
\label{eq:bound-h-0}
&\displaystyle M_1=\sup_{h\in\mathscr H}h(0)<\infty,& \\
\label{eq:unif-int}
&\displaystyle \lim_{x\to\infty}\sup_{h\in\mathscr H } \int_x^\infty
h(y)\,dy=0.&
\end{eqnarray}
Let $U$ be the renewal function associated with a nonlattice
inter-renewal distribution with finite mean $\beta$. Then
%
%
\begin{equation}
\label{eq:2}
\lim_{x\to\infty}\sup_{h\in\mathscr H}
\biggl|U*h(x)- \frac{1}{\beta} \int_0^\infty h(y)\,dy\biggr|=0.
\end{equation}
\end{lem}
\begin{pf}
Conditions (\ref{eq:bound-h-0}) and (\ref{eq:unif-int}) imply that
\[
M_2\equiv\sup_{h\in\mathscr H } \int_0^\infty h(y)\,dy<\infty.
\]
Let $\delta>0$ and $\varepsilon>0$ be arbitrary positive numbers in
$(0, 1)$. By (\ref{eq:unif-int}), there exists $N=N(\delta, \varepsilon
)$ such
that for each $h\in\mathscr H$
\[
\int_{N\delta}^\infty h(y)\,dy < \delta\varepsilon.
\]
Furthermore, by the Blackwell theorem, there exists $x^*$ such that
for each $x\ge x^*$ and each $k=0, \ldots, N$,
\[
\frac{\delta}{\beta}-\delta\varepsilon
< U\bigl(x-(k+1)\delta\bigr) - U(x-k\delta)
< \frac{\delta}{\beta} + \delta\varepsilon.
\]
Let $d_k(x)=U(x-k\delta)-U(x-(k+1)\delta)$ for all $k\ge0$. Here,
we take the convention that $U(x)=0$ for all $x<0$. Define
\[
\bar h^{\delta}(x) = \sum_{k=0}^\infty h(k\delta)
1_{\{k\delta\le x < (k+1)\delta\}}.
\]
Clearly, $h(x)\le\bar h^{\delta}(x)$ for $x\ge0$. So for
$x\ge x^*$ and $h\in\mathscr H$, we have
\begin{eqnarray*}
U*h(x)&\le& U*{\bar{h}^\delta}(x) =\sum_{k=0}^\infty
h(k\delta)d_k(t)\\
& =& \sum_{k=0}^N h(k\delta)d_k(t) + \sum_{k=N+1}^\infty
h(k\delta)d_k(t) \\
& \le& \sum_{k=0}^N h(k\delta)\biggl(\frac{\delta}{\beta}+\varepsilon
\delta\biggr) +
U(\delta) \sum_{k=N+1}^\infty
h(k\delta) \\
& \le&\biggl(\delta h(0)+\int_0^\infty h(y)\,dy\biggr)\biggl(\frac{1}{\beta}+\varepsilon\biggr) +
U(\delta) \frac{1}{\delta} \int_{N\delta}^\infty h(y)\,dy \\
& \le&\biggl(\delta h(0)+\int_0^\infty h(y)\,dy\biggr)\biggl(\frac{1}{\beta}+\varepsilon\biggr) +
U(1) \varepsilon\\
& \le& \frac{1}{\beta}\int_0^\infty h(y)\,dy
+ \delta M_1\biggl(\frac{1}{\beta}+\varepsilon\biggr)+
M_2\varepsilon+ U(1) \varepsilon.
\end{eqnarray*}
Define
\[
\underline h^{\delta}(t)
= \sum_{k=0}^\infty h\bigl((k+1)\delta\bigr)1_{\{k\delta\le t < (k+1)\delta\}}.
\]
Clearly, $h(x)\ge\underline h^{\delta}(x)$ for $x\ge0$. So for
$x\ge x^*$ and $h\in\mathscr H$, we have
\begin{eqnarray*}
U*h(x)&\ge& U*{\underline{h}^\delta}(x) =\sum_{k=0}^\infty
h\bigl((k+1)\delta\bigr)d_k(t)\\
&=& \sum_{k=0}^N h\bigl((k+1)\delta\bigr)d_k(t) + \sum_{k=N+1}^\infty
h\bigl((k+1)\delta\bigr)d_k(t) \\
&\ge&\sum_{k=0}^N h\bigl((k+1)\delta\bigr)\biggl(\frac{\delta}{\beta}-\varepsilon
\delta
\biggr)\\
&=& \Biggl(
\sum_{k=0}^\infty h(k\delta) -\delta h(0) -
\sum_{k=N+2}^\infty h(k\delta)
\Biggr)\biggl(\frac{\delta}{\beta}-\varepsilon
\delta\biggr)\\
& \ge&\biggl(\int_0^\infty h(y)\,dy -
\int_{(N+1)\delta}^\infty h(y)\,dy \biggr)\biggl(\frac{1}{\beta}-\varepsilon\biggr)
- \delta h(0) \biggl(\frac{\delta}{\beta}-\varepsilon\delta\biggr) \\
&\ge& \frac{1}{\beta}\int_0^\infty h(y)\,dy
- M_2\varepsilon
- \delta\varepsilon\biggl(\frac{1}{\beta}-\varepsilon\biggr)
- \delta^2 M_1 \biggl(\frac{1}{\beta}-\varepsilon\biggr).
\end{eqnarray*}
Thus,
\begin{eqnarray*}
&& \limsup_{x\to\infty}\sup_{h\in\mathscr H} \biggl| U*h(x) -
\frac{1}{\beta}\int_0^\infty
h(y)\,dy \biggr| \\
&&\qquad \le
\delta M_1\biggl(\frac{1}{\beta}+\varepsilon\biggr)
+2M_2\varepsilon+ U(1) \varepsilon
+\delta\varepsilon\biggl(\frac{1}{\beta}-\varepsilon\biggr)\\
&&\qquad\quad{}+\delta^2 M_1 \biggl(\frac{1}{\beta}-\varepsilon\biggr).
\end{eqnarray*}
Because $\delta>0$ and $\varepsilon>0$ can be arbitrarily small, we
have
\[
\lim_{x\to\infty}\sup_{h\in\mathscr H}
\biggl| U*h(x) - \frac{1}{\beta}\int_0^\infty h(y)\,dy \biggr|=
0.
\]
\upqed\end{pf}

\section{Some results on the Prohorov metric}
\label{sec:conv-proh-metr}

Lemma \ref{lem:conv-proh} is applied in
Section \ref{subsec:unif-conv-invar-manif}, and
Lemma \ref{lem:workload-pov} is applied in
Section \ref{subsec:proof-ssc}. Since we could not find these results
in the literature, we include them here for completeness.
\begin{lem}\label{lem:conv-proh}
Let $\mu$ and $\mu_1$ be finite Borel measures on $[0,\infty)$.
Denote $A_y=(y,\infty)$ for all $y\ge0$. Let
$M=\max(\langle{\chi},{\mu}\rangle,\langle{\chi},{\mu_1}\rangle
)$. For all $0<\varepsilon<1$
if
%
%
\begin{equation}\label{eq:proh-cond}
\sup_{y\ge0}|\mu(A_y)-\mu_1(A_y)|<\varepsilon,
\end{equation}
then
\[
\pov[\mu,\mu_1]<(M+2)\varepsilon^{1/3}.
\]
\end{lem}
\begin{pf}
Let $\alpha, \beta$ be positive constants to be determined later.
Note that
\[
\mu((\varepsilon^{-\alpha},\infty))\le M\varepsilon^\alpha.
\]
For any real number $a$, denote $I_a=(a,a+\varepsilon^\beta]$.
Condition (\ref{eq:proh-cond}) implies that
\[
\sup_{a\in\R}|\mu(I_a)-\mu_1(I_a)|<2\varepsilon.
\]
For any Borel set $A\subset[0,\infty)$, there exist $a_1,\ldots,
a_N$ such that
\[
A\cap[0,\varepsilon^{-\alpha}]\subset\bigcup_{i=1}^N I_{a_i},
\]
and $I_{a_i}\cap I_{a_j}=\varnothing$ for all $i\ne j$, and
$I_{a_i}\cap A\ne\varnothing$ for all $i$. These conditions imply
that
\[
N\le\varepsilon^{-\alpha-\beta}
\]
and
\[
\bigcup_{i=1}^N I_{a_i}\subset A^{\varepsilon^\beta},
\]
where $A^{\varepsilon^\beta}$ is the $\varepsilon^\beta$-enlargement of
the set defined in Section \ref{subsec:notatioin}. So we have
\begin{eqnarray*}
\mu(A)&\le&\mu(A\cap[0,\varepsilon^{-\alpha}])
+\mu\bigl(A\cap(\varepsilon^{-\alpha},\infty)\bigr)\\
&\le&\mu\Biggl(\bigcup_{i=1}^N I_{a_i}\Biggr)+M\varepsilon^\alpha\\
&<&\mu_1\Biggl(\bigcup_{i=1}^N I_{a_i}\Biggr)+N2\varepsilon+M\varepsilon^\alpha\\
&\le&\mu_1(A^{\varepsilon^\beta})+2\varepsilon^{1-\alpha-\beta
}+M\varepsilon
^\alpha.
\end{eqnarray*}
Now choose $\alpha=\beta=1/3$ to obtain
\[
\mu(A)<\mu_1\bigl(A^{(M+2)\varepsilon^{1/3}}\bigr)+(M+2)\varepsilon^{1/3}.
\]
Exchanging the position of $\mu$ and $\mu_1$ in the above argument,
we have
\[
\mu_1(A)<\mu\bigl(A^{(M+2)\varepsilon^{1/3}}\bigr)+(M+2)\varepsilon^{1/3}.
\]
This completes the proof.
\end{pf}
\begin{lem}\label{lem:workload-pov}
Suppose $\mu_1$ and $\mu$ are finite Borel measures on $\R_+$
satisfying
%
%
\begin{equation}\label{eq:tech-6-cond-1}
\pov[\mu_1,\mu] <\varepsilon<1,
\end{equation}
and $\langle{\chi^{1+q}},{\mu_1}\rangle<M$, $\langle{\chi
^{1+q}},{\mu}\rangle<M$ for some
positive constants $q$ and $M$, then
\[
|\langle{\chi},{\mu_1}\rangle-\langle{\chi},{\mu}\rangle|\le
\varepsilon^{1/2}+\frac
{2M}{q}\varepsilon^{q/2}.
\]
\end{lem}
\begin{pf}
By Markov inequality, $\mu_1(A_x)\le\frac{M}{x^{1+q}}$ and
$\mu(A_x)\le\frac{M}{x^{1+q}}$ for all \mbox{$x\ge0$}. For any $C>0$, we
have the following inequality:
\begin{eqnarray*}
|\langle{\chi},{\mu_1}\rangle-\langle{\chi},{\mu}\rangle|
&\le&\int_0^C|\mu_1(A_x)-\mu(A_x)|\,dx\\
&&{}+\int_C^\infty\mu_1(A_x)\,dx+\int_C^\infty\mu(A_x)\,dx\\
&\le& C\varepsilon+2\int_C^\infty\frac{M}{x^{1+q}}\,dx\\
&=&C\varepsilon+2M\frac{1}{qC^{q}}.
\end{eqnarray*}
The result follows by letting $C=\varepsilon^{-{1/2}}$.
\end{pf}

\vspace*{-14pt}

\section{Glivenko--Cantelli estimate}
\label{append:glivenko-cantelli}

For any $r$, consider the sequence of i.i.d. random variables
$\{v^r_i\}_{i=-\infty}^{\infty}$ with law~$\nu^r$. In our setting,
those $v^r_i$'s with $i\ge1$ correspond to the service requirement of
the arriving jobs in the $r$th system; those with $i\le0$ correspond
to the service requirement of initial jobs waiting in the buffer. For
any $n\in\Z$ and $l\in\R_+$, define
%
%
\begin{equation}\label{eq:def-eta}
\fr\eta(n,l)=\frac{1}{r}\sum_{i=n+1}^{n+\lfloor{rl}\rfloor}\delta
_{v^r_i}.
\end{equation}
The objective of this section is to obtain the \textit{Glivenko--Cantelli
estimate}, Lem\-ma~\ref{lem:glivenko-cantelli} below, for
$\fr\eta(n,l)$. Very similar result was shown in Lemma 4.7
\cite{GromollKurk2007}. For completeness, the proof which follows the
one in \cite{GromollKurk2007} is provided here.

To present the result, we introduce some notions from empirical
process theory. Our primary references are \cite{GromollKurk2007} and
\cite{Varrt1996}.

A collection $\testset$ of subsets of $\R^2$ shatters an $n$-point
subset $\{x_1,\ldots,x_n\}\subset\R_+$ if the collection
$\{\testset\cap\{x_1,\ldots,x_n\}\dvtx C\in\testset\}$ has cardinality
$2^n$. In this case, we say that $\testset$ picks out all subsets of
$\{x_1,\ldots,x_n\}$. The \textit{Vapnik--\v{C}ervonenkis index}
(\textit{VC-index}) of $\testset$ is
\[
V_\testset=\min\{n\dvtx\testset\mbox{ shatters no $n$-point
subset}\},
\]
where the minimum of the empty set equals infinity.
The collection $\testset$ is a \textit{Vapnik--\v{C}ervonenkis class}
(\textit{VC-class}) if it has finite VC-index. Let $\testfn$ be a family of
Borel measurable functions $f\dvtx R_+\to\R$. We call $\testfn$ a VC-class
if the collection of subgraphs $\{\{(x,y)\dvtx y<f(x)\}\dvtx f\in\testfn
\}$ is a
VC-class of sets in $\R^2$.

We call a family of functions $\testfn$ a \textit{Borel measurable
class} if, for each $n\in\mathbb N$ and
$(e_1,\ldots,e_n)\in\{-1,1\}^n$, the map
\[
(x_1,\ldots,x_n)\to\sup_{f\in\testfn}\sum_{i=1}^{n}e_if(x_i)
\]
is Borel measurable on $\R_+^n$. The condition requires that, for all
$\delta>0$ and $r\in\R_+$, the families
$\testfn_\delta^r=\{f-g\dvtx f,g\in\testfn,\|f-g\|_{\nu^r,2}<\delta\}$ and
$\testfn_\infty^2=\{(f-g)^2\dvtx f,g\in\testfn\}$ are Borel measurable,
where
\[
\|f\|_{\nu^r,2}=\langle{|f|^2},{\nu}\rangle^{1/2}
\]
denotes the $L_2(\nu^r)$-norm.

We call a Borel measurable function $\bar f\dvtx R_+\to\R$ an envelope
function for $\testfn$ if any element in $\testfn$ is bounded by
$\bar
f$. A VC-class with an envelop function satisfies a very useful
entropy bound. Let $\Gamma$ be the set of finitely discrete
probability measures $\gamma$ on $\R_+$ such that $\|\bar
f\|_{\gamma,2}>0$. For any Borel measurable function $f\dvtx\R_+\to\R$
satisfying $\|f\|_{\gamma,2}<\infty$, let
$B_f(\varepsilon)=\{g\in\testfn\dvtx\|g-f\|_{\nu^r,2}<\varepsilon\}$
denote the
$L_2(\nu^r)$-ball in $\testfn$, centered at $f$ with radius
$\varepsilon$. For a family of functions $\testfn$,
$N(\varepsilon,\testfn,L_2(\gamma))$ is the smallest number of balls
$B_f(\varepsilon)$ needed to cover $\testfn$. Then $\testfn$ satisfies
%
%
\begin{equation}\label{ineq:entropybound}
\int_0^\infty\sup_{\gamma\in\Gamma}
\sqrt{\log N(\varepsilon\|\bar f\|_{\gamma,2},\testfn,L_2(\gamma))}\,
d\varepsilon<\infty;
\end{equation}
see Definition 2.1.5, (2.5.1) and Theorem 2.6.7 in
\cite{Varrt1996}.
%
%
\begin{lem}\label{lem:glivenko-cantelli}
Let $\testfn$ be a VC-class of Borel measurable functions such that
$\testfn^2_\infty$ and $\testfn^r_\delta$ are Borel measurable
classes for all $r\in\R_+$ and $\delta>0$. Assume there exists an
envelop function $\bar f$ of $\testfn$ such that
%
%
\begin{equation}\label{ineq:envelop}
\lim_{N\to\infty}\sup_{r\in\R_+}\langle{\bar f^21_{\bar
f>N}},{\nu^r}\rangle=0.
\end{equation}
Fix constants $M_1,L_1>0$. For all $\varepsilon,\varepsilon'>0$,
%
%
\begin{eqnarray}\label{ineq:glivenko-cantelli}
&&\limsup_{r\to\infty}\mathbb{P}^r\Bigl({\max_{-rM_1< n< r^2M_1}\sup
_{l\in[0,L_1]}\sup_{f\in\testfn} |\langle{f},{\fr\eta
(n,l)}\rangle-l\langle{f},{\nu^r(A_x)}\rangle|>\varepsilon' }
\Bigr)\hspace*{-28pt}\nonumber\\[-8pt]\\[-8pt]
&&\qquad<\varepsilon.\nonumber
\end{eqnarray}
\end{lem}
%
%
\begin{rem}
To apply the lemma in this paper, we take
%
%
\begin{equation}\label{eq:set-of-test-fns}
\testfn=\{1_C\dvtx C\in\testset\}\cup\{\chi^{1+p},\chi^{2+p}\},
\end{equation}
where $\testset=\{[y,\infty)\dvtx y\in\R_+\}\cup\{(y,\infty)\dvtx y\in\R
_+\}$
and $p$ is the same one as in condition (\ref{eq:cond-S-u}). It is
very easy to see that $\testfn$ is a VC-class. Note that both
$\testfn_\delta^r$ and $\testfn_\infty^r$ are subsets of functions
of the form $1_{(a,b)}$ (or $1_{(a,b]}$, $1_{[a,b)}$, $1_{[a,b]}$).
So the supreme over these two families will be the same as supreme
over all $a,b$ in subsets of $\R_+$, which will be the same as over
all $a,b$ in subsets of $Q$. Borel measurability is preserved when
take supreme over a countable set. It is also clear that
\[
\bar f(x)=\cases{
1, &\quad $x<1$,\cr
x^{2+p}, &\quad $x\ge1$,}
\]
is an envelop function. Condition (\ref{eq:cond-S-u}) implies
(\ref{ineq:envelop}).
\end{rem}

To better structure the proof, we present the following auxiliary
lemma.
%
%
\begin{lem}\label{lem:empirical}
For $n\in\Z$ and $k\in\mathbb N$, define
%
%
\begin{equation}\label{eq:def-of-xi-nk}
\xi^r_{n,k}=\frac{1}{\sqrt k}\sum_{i=n+1}^{n+k}(\delta_{v^r_i}-\nu^r).
\end{equation}
%
Then for any $q>1$, $y>2$ and $n\in\Z$ there exists $M_q<\infty$ and
$k_0$ such that $k\ge k_0$ implies
%
%
\begin{equation}\label{ineq:boundempirical}
\sup_{r}\mathbb{P}^r\Bigl({\sup_{f\in\testfn}\langle{f},{\xi
^r_{n,k}}\rangle>y}\Bigr)<\frac{M_q}{y^q}.
\end{equation}
The constant $M_q$ does not depend on $y$.
\end{lem}
\begin{pf}
Let us first fix $n=0$ and look at $\xi^r_{0,k}$ which will be
denoted by $\xi^r_k$ for simplicity. The property
(\ref{ineq:envelop}) of the envelop function $\bar f$ and the
uniform entropy bound (\ref{ineq:entropybound}), together with the
sets $\testfn_\delta^r$ and $\testfn_\infty^r$ being Borel
measurable, imply that $\testfn$ is Donsker and pre-Gaussian
uniformly in $\nu^r$, $r\in\R_+$. (See Theorem 2.8.3 in~\cite{Varrt1996}.)

Let $l^\infty(\testfn)$ be the space of all probability measures on
$\R_+$ equipped with norm
$\|\cdot\|_{\testfn}=\sup_{f\in\testfn}\langle{f},{\cdot}\rangle
$. $\testfn$
being Donsker uniformly in $\nu^r$ means that $\xi^r_k$ converges
weakly as $n\to\infty$ in $l^\infty(\testfn)$ to a tight, Borel
measurable version of the Brownian bridge $\xi^r$ uniformly for all
$\nu^r$. According to Chapter $1.12$ in \cite{Varrt1996}, this is
equivalent to
%
%
\begin{equation}\label{eq:BL1}
{\sup_{h\in\mathrm{BL}_1}}|\E^rh(\xi^r_k)-\E h(\xi^r)|\to0.
\end{equation}
uniformly for all $\nu^r$, where $\mathrm{BL}_1$ is the set of
functions $h\dvtx l^\infty(\testfn)\to\R$ which are uniformly bounded by
$1$ and satisfy $|h(z_1)-h(z_2)|\le\|z_1-z_2\|_{\testfn}$.
Pre-Gaussian uniformly in $\nu^r$ means that
%
%
\begin{equation}\label{ineq:preGaussian}
\sup_{r}\E^r\Bigl[\sup_{f\in\testfn}\langle{f},{\xi^r}\rangle
\Bigr]<\infty.
\end{equation}

Define $h_y\dvtx l^\infty(\testfn)\to\R$ by
\[
h_y(\cdot)=\Bigl(\sup_{f\in\testfn}\langle{f},{\cdot}\rangle
-y+1\Bigr)^+\land1.
\]
Then it is clear that $h_y\in\mathrm{BL}_1$, and
\[
\sup_{r}\mathbb{P}^r\Bigl({\sup_{f\in\testfn}\langle{f},{\xi
^r_k}\rangle>y}\Bigr)
\le\sup_{r}\E^r[h_y(\xi^r_k)].
\]
By (\ref{eq:BL1}) and the above inequality, there exists
$k_0\in\mathbb N$ such that $k\ge k_0$ implies
\[
\sup_{r}\mathbb{P}^r\Bigl({\sup_{f\in\testfn}\langle{f},{\xi
^r_k}\rangle>y}\Bigr)
\le\sup_{r}\E^r[h_y(\xi^r)]+y^{-q}.
\]
Applying the definition of $h_y$ and Markov inequality to obtain
\begin{eqnarray*}
\sup_{r}\mathbb{P}^r\Bigl({\sup_{f\in\testfn}\langle{f},{\xi
^r_k}\rangle>y}\Bigr)
&\le&\sup_{r}\mathbb{P}^r\Bigl({\sup_{f\in\testfn}\langle
{f},{\xi^r}\rangle>y-1}\Bigr)+y^{-q}\\
&\le& y^{-q}
\Bigl(2^q\sup_{r}\E^r\Bigl[\sup_{f\in\testfn}\langle{f},{\xi
^r}\rangle
\Bigr]^q+1\Bigr).
\end{eqnarray*}
Let $M_q$ be the last term in parentheses, which does not depend on
$y$. For each $r\in\R_+$, the Brownian bridge is separable and
Gaussian with $\sup_{f\in\testfn}\langle{f},{\xi^r}\rangle$
finite almost
surely. Thus, there exist a constant $M$ such that for all
$r\in\R_+$,
\[
\E^r\Bigl[\sup_{f\in\testfn}\langle{f},{\xi^r}\rangle\Bigr]^q
\le M\Bigl[\E^r\sup_{f\in\testfn}\langle{f},{\xi^r}\rangle\Bigr]^q,
\]
see Proposition A.2.4 in \cite{Varrt1996}. Conclude from
(\ref{ineq:preGaussian}) that $M_q<\infty$.

So far, we have shown that the result (\ref{ineq:boundempirical}) is
true for $n=0$. Note that for any $n\in\Z$, $\xi^r_{n,k}$ is
defined on the shifted sequence $v_{n+1}^r,v_{n+2}^r,\ldots.$ By the
i.i.d. property of the sequence, if we fix $k$ then $\xi^r_{n,k}$ has
the same distribution for all $n\in\Z$. So we can conclude that
(\ref{ineq:boundempirical}) is true for all $n\in\Z$.
\end{pf}
%
%
\begin{pf*}{Proof of Lemma \ref{lem:glivenko-cantelli}}
Note that
\[
|\langle{f},{\fr\eta(n,l)}\rangle-l\langle{f},{\nu^r}\rangle|\le
\frac{1}{r}\sum_{i=n+1}^{n+\lfloor{rl}\rfloor}
[\langle{f},{\delta_{v^r_i}}\rangle-\langle{f},{\nu^r}\rangle
]+\frac{1}{r}.
\]
Since for each $\varepsilon'>0$, $1/r<\varepsilon'/2$ for all large $r$,
so the probability in (\ref{ineq:glivenko-cantelli}) can be bounded
by
%
%
\begin{equation}\label{eq:tech-A-1}
\limsup_{r\to\infty}\mathbb{P}^r\Biggl({\max_{-rM_1< n< r^2M_1}\sup
_{l\in[0,L]}\sup_{f\in\testfn} \Biggl|\frac{1}{r}\sum
_{i=n+1}^{n+\lfloor{rl}\rfloor}[\langle{f},{\delta_{v^r_i}}\rangle
-\langle{f},{\nu^r}\rangle]\Biggr|>\frac{\varepsilon
'}{2} }\Biggr).\hspace*{-38pt}
\end{equation}
Pick $\delta>0$, when $r$ is large enough ($r>M_1/\delta$) the
interval $[-rM_1,r^2M_1]$ will be covered by intervals
\[
[-r^2\delta, 0], [0,r^2\delta], \ldots,
\biggl[\biggl(\biggl\lceil{\frac{M_1}{\delta}}\biggr\rceil-1\biggr)r^2\delta,\biggl\lceil{\frac
{M_1}{\delta}}\biggr\rceil r^2\delta\biggr].
\]
When $r$ is large enough ($r^2\delta>\lfloor{rL_1}\rfloor$),
(\ref{eq:tech-A-1}) can be further bounded by
\begin{eqnarray*}
&&\limsup_{r\to\infty}\mathbb{P}^r\biggl(\max_{-1\le j \le\lceil
{{M_1}/{\delta}}\rceil-1}\max_{0\le k,k'\le r^2\delta} \sup
_{f\in
\testfn}\biggl|\frac{\sqrt{k}}{r}\langle{f},{\xi^r_{jr^2\delta
,k}}\rangle\\
&&\qquad\hspace*{160.5pt}{}-\frac
{\sqrt{k'}}{r}\langle{f},{\xi^r_{jr^2\delta,k'}}\rangle\biggr| >\frac
{\varepsilon
'}{2} \biggr).
\end{eqnarray*}
Since $\xi^r_{n,\cdot}$ has stationary increments, the previous term
can be bounded above by
\[
\limsup_{r\to\infty}\biggl\lceil{\frac{M_1}{\delta}}\biggr\rceil\mathbb
{P}^r
\biggl({\max_{0\le k\le r^2\delta}\sup_{f\in\testfn} \biggl|\frac{\sqrt
{k}}{r}\langle{f},{\xi^r_{0,k}}\rangle\biggr|>\frac{\varepsilon'}{2}}\biggr).
\]
By Ottaviani's inequality (see Proposition A.1.1 in
\cite{Varrt1996}) and by stationary increments of $\xi^r_{n,\cdot}$,
this can be bounded above by
%
%
\begin{equation}\label{eq:ottaviani}
\limsup_{r\to\infty}\frac{\lceil{{M_1}/{\delta}}\rceil
\mathbb{P}^r({\sup_{f\in\testfn}\langle{f},{\xi
^r_{0,\lfloor{r^2\delta}\rfloor}}\rangle>{\varepsilon
'}/({4\sqrt\delta}) })
}
{1-\max_{0\le k\le r^2\delta}\mathbb{P}^r({\sup_{f\in\testfn
}\langle{f},{\xi^r_{0,k}}\rangle>{\varepsilon'r}/({4\sqrt k}) })
}.
\end{equation}
Assume $\delta$ is small enough so that
$\frac{\varepsilon'}{4\sqrt\delta}>2$. By Lemma \ref{lem:empirical},
there exists $M_3$ and $k_0\in\N$ such that $k>k_0$ implies
\[
\sup_{r}\mathbb{P}^r\biggl({\sup_{f\in\testfn}\langle{f},{\xi
^r_{0,k}}\rangle>\frac{\varepsilon'}{4\sqrt\delta} }\biggr)\le
\biggl(\frac
{4\sqrt\delta}{\varepsilon'}\biggr)^3M_3.
\]
Since $\lfloor{r^2\delta}\rfloor\to\infty$ as $r\to\infty$, the
limit superior
of the numerator in (\ref{eq:ottaviani}) can be bounded above by
$\lceil{M_1/\delta}\rceil({4\sqrt\delta}/{\varepsilon'})^3M_3$,
which can be
made arbitrarily small by choosing $\delta$ sufficiently small. By
the same reason, those terms in the maximum of the denominator with
index $k>k_0$ are bounded above by
$({4\sqrt\delta}/{\varepsilon'})^3M_3$. For those terms with index
$k\le k_0$,
\[
\mathbb{P}^r\biggl({\sup_{f\in\testfn}\langle{f},{\xi
^r_{0,k}}\rangle>\frac
{\varepsilon'r}{4\sqrt k} }\biggr)
\le\mathbb{P}^r\biggl({\sup_{f\in\testfn}\langle{f},{\xi
^r_{0,k}}\rangle>\frac{\varepsilon'r}{4\sqrt k_0} }\biggr),
\]
which converges to zero as $r\to\infty$. By choosing $\delta$ small
enough, (\ref{eq:ottaviani}) can be made arbitrarily small for all
large $r$.
\end{pf*}
\end{appendix}

\section*{Acknowledgments}

The authors thank two anonymous referees for significantly improving
the paper.


%
\printaddresses

\end{document}